\documentclass[twoside,reqno,a4paper]{amsart}

\usepackage{style}

\begin{document}
\title[de Branges systems]{De {B}ranges canonical systems with finite \\logarithmic integral}
\author{Roman V.~Bessonov, \quad  Sergey A.~Denisov}

\address{
\begin{flushleft}
Roman Bessonov: bessonov@pdmi.ras.ru\\\vspace{0.1cm}
St.\,Petersburg State University\\  
Universitetskaya nab. 7-9, 199034 St.\,Petersburg, RUSSIA\\
\vspace{0.1cm}
St.\,Petersburg Department of Steklov Mathematical Institute\\ Russian Academy of Sciences\\
Fontanka 27, 191023 St.Petersburg,  RUSSIA\\
\end{flushleft}
}
\address{
\begin{flushleft}
Sergey Denisov: denissov@wisc.edu\\\vspace{0.1cm}
University of Wisconsin--Madison\\  Department of Mathematics\\
480 Lincoln Dr., Madison, WI, 53706, USA\vspace{0.1cm}\\
\vspace{0.1cm}
Keldysh Institute of Applied Mathematics\\ Russian Academy of Sciences\\
Miusskaya pl. 4, 125047 Moscow, RUSSIA\\
\end{flushleft}
}

\thanks{
The work of RB in Sections 5 and 7 is supported by grant RScF 19-71-30002 of the Russian Science Foundation.
The work of SD done in Sections 9 and 10 is supported by grants RScF-14-21-00025 and RScF-19-71-30004 of the Russian Science Foundation. His research conducted in the rest of the paper is supported by the grants NSF-DMS-1464479,  NSF DMS-1764245, and Van Vleck Professorship Research Award. SD gratefully acknowledges the hospitality of IHES where part of this work was done.}
\subjclass[2010]{42C05, 34L40, 34A55}
\keywords{Szeg\H{o} class, Canonical Hamiltonian system, Inverse problem, Entropy.}

\begin{abstract} Krein -- de Branges spectral theory establishes a correspondence between the class of differential operators called canonical Hamiltonian systems and measures on the real line with finite Poisson integral. We further develop this area by giving a description of canonical Hamiltonian systems whose spectral measures have  logarithmic integral converging over the real line. This result can be viewed as a spectral version of the classical Szeg\H{o} theorem in the theory of polynomials orthogonal on the unit circle. It extends Krein--Wiener completeness theorem, a key fact in the prediction of stationary Gaussian processes.    
\end{abstract}

\maketitle

\setcounter{tocdepth}{1}

\tableofcontents


\pagebreak

\section{Introduction}\label{s1}
In this paper, we look at the spectral theory of de Branges' canonical system, which is defined by the system of differential equations of the form\smallskip
\begin{equation}\label{eq1}
J \tfrac{d}{dt} M(t,z) = z \Hh(t) M(t,z), \quad M (0, z) = I_{2\times 2}\dd \idm, \quad J\dd \jm,\,\quad t \ge 0,\quad  z \in \C\,.
\end{equation}\smallskip
The $2\times 2$ matrix-function $\Hh$ on $\R_+= [0,+\infty)$ is called the Hamiltonian of canonical system \eqref{eq1}. We will always assume that $\Hh$ satisfies the following conditions:
\begin{enumerate}
\item[$(a)$]$\Hh(t) \ge 0$ and \,$\trace \Hh(t) > 0$ for Lebesgue almost every $t \in \R_+$,
\item[$(b)$] the entries of $\Hh$ are real measurable functions absolutely integrable on compact subsets of $\R_+$.
\end{enumerate}
In 1960's, L. de Branges developed his theory of Hilbert spaces of entire functions (see \cite{dbbook} and \cite{Romanov,Remlingb} for recent exposition). One result of this monumental work is the theorem that establishes a bijection between Hamiltonians $\Hh$ in \eqref{eq1} and nonconstant analytic functions in $\C_+ = \{z \in \C:\; \Im z> 0\}$ with nonnegative imaginary part. Every such  function is generated by a nonnegative measure on the real line. In this paper, we make a further step in de Branges' theory by identifying Hamiltonians that correspond to measures in the Szeg\H{o} class, i.e., the measures whose logarithmic integral converges over~$\R$.

\medskip

To formulate the main results of the paper, we need some definitions. A Hamiltonian $\Hh$ on $\R_+$ is called singular if 
$$
\int_{0}^{+\infty}\trace\Hh(t)\,dt = +\infty.
$$ 
Two  Hamiltonians $\Hh_1$, $\Hh_2$ on $\R_+$ are called equivalent if there
exists an increasing absolutely continuous function $\eta$ defined on $\R_+$ such that $\eta(0) = 0$, $\lim_{t \to+\infty}\eta(t) = +\infty$, and $\Hh_2(t) = \eta'(t)\Hh_1(\eta(t))$ for Lebesgue almost every $t \in \R_+$. Clearly, $\eta(t)$ rescales the variable $t$. We say
that Hamiltonian $\Hh$ is trivial if there is a  non-negative matrix $A$ with $\rank A=1$,  such that $\Hh$ is equivalent to $A$, i.e., $\Hh(t) =\eta'(t) A$ for a.e. $t\in \R_+$, where $\eta$ is an increasing absolutely continuous function on $\R_+$, which satisfies $\eta(0) = 0$ and $\lim_{t \to+\infty}\eta(t) = +\infty$. If Hamiltonian is not trivial, it is called nontrivial.

\medskip
 
Recall that function $m$ belongs to the Herglotz-Nevanlinna class $\mathcal{N}(\C_+)$   if it is analytic in $\C_+$ and $\Im m(z)\ge 0$ for $z\in \C_+$. It is well-known \cite{Garnett}, that $m\in \mathcal{N}(\C_+)$ if and only if it admits the following  representation
\begin{equation}\label{hg}
m(z) = \frac{1}{\pi}\int_{\R}\left(\frac{1}{x-z} - \frac{x}{x^2+1}\right)\,d\mu(x)+ bz+a, \quad z\in \mathbb{C}_+,
\end{equation}
where $b \ge 0$, $a \in\R$, and $\mu$ is a Radon measure on $\R$, which satisfies
\begin{equation}\label{eq111}
\int_{\R}\frac{d\mu}{1+x^2}<\infty.
\end{equation}
We call measures on $\R$ satisfying \eqref{eq111} Poisson-finite. The class  $\mathcal{N}(\C_+)$  appears naturally in the theory of canonical Hamiltonian systems. Let $\Hh$ be a nontrivial and singular Hamiltonian.  Given condition $(b)$ on $\Hh$, there exists unique matrix-valued function  $M$ that solves \eqref{eq1}. Denote by $\Theta^\pm$, $\Phi^{\pm}$ its entries so that
\begin{equation}\label{eq32}
M(t,z) = (\Theta(t,z), \Phi(t,z)) = \begin{pmatrix}\Theta^{+}(t,z) & \Phi^{+}(t,z) \\ \Theta^{-}(t,z) &\Phi^{-}(t,z)\end{pmatrix}.
\end{equation}
Fix a parameter $\omega \in \R \cup \{\infty\}$. The Titchmarsh-Weyl function of $\Hh$ is defined by
\begin{equation}\label{mf}
m(z) = \lim_{t \to +\infty}\frac{\omega\Phi^{+}(t,z) + \Phi^{-}(t,z)}{\omega\Theta^{+}(t,z) + \Theta^{-}(t,z)}, \qquad z \in \C_+,
\end{equation}
where the fraction $\frac{\infty c_1 + c_2}{\infty c_3 + c_4}$ for non-zero numbers $c_1$, $c_3$ is interpreted as $c_1/c_3$. In Weyl's theory for canonical systems (see \cite{HSW} or
Section 8 in \cite{Romanov}), it is shown  that the expression under the
limit in \eqref{mf} is well-defined for large $t > 0$ (i.e., the denominator is non-zero) for every given singular nontrivial Hamiltonian $\Hh$. Moreover, the limit $m(z)$ exists, does not depend on $\omega$, $m$ is analytic in $z\in \C_+$ and has non-negative imaginary part, i.e., $m\in \cal{N}(\C_+)$. In particular, $m$ admits   representation \eqref{hg}.  The measure $\mu$ in \eqref{hg} is called the spectral measure for the Hamiltonian $\Hh$. It is easy to check that equivalent Hamiltonians have equal Titchmarsh-Weyl functions, see \cite{WW12}.

\medskip

 Now we can formulate the result of de Branges that establishes a bijection between Hamiltonians and Herglotz-Nevanlinna functions. See  \cite{dbbook},  \cite{Romanov}, \cite{Winkler95} and also \cite{KK68}  for its proofs.

\begin{Thm}{\rm (de Branges)}
For every nonconstant function $m \in \mathcal{N}(\C_+)$, there exists a singular nontrivial Hamiltonian $\Hh$ on $\R_+$ such that $m$ is the Titchmarsh-Weyl function \eqref{mf} for $\Hh$. Moreover, any two singular nontrivial Hamiltonians $\Hh_1$, $\Hh_2$ on $\R_+$ generated by $m$ are equivalent.
\end{Thm}
For trivial Hamiltonians, function $m$ is a real constant. Indeed, in that case, one can solve \eqref{eq1} explicitly and this calculation shows that $m(z)=const\in \R\cup\infty$. For example, $\Hh=\left(\begin{smallmatrix} 1&0\\0&0\end{smallmatrix}\right)$ gives 
\begin{equation}\label{sd_ff90}
\Theta^+=1,\quad\Theta^-=-zt,\quad\Phi^+=0,\quad \Phi^-=1,
\end{equation}
 so $m=0$. Similarly, if $\Hh=\left(\begin{smallmatrix} 0&0\\0&1\end{smallmatrix}\right)$, then $\Theta^+=1,\Theta^-=0,\Phi^+=zt,\Phi^-=0$ and we  let $m=\infty$. 

\medskip

Given a Poisson-finite measure $\mu$ on $\R$, we will denote by $w$ the
density of $\mu$ with respect to the Lebesgue measure $dx$ on $\R$, and by $\mus$  the singular part of $\mu$, so
that $\mu = w\,dx+\mus$.
In this paper, our aim  is to characterize singular nontrivial Hamiltonians whose spectral measures
have finite logarithmic integral, i.e., the integral
$$
\int_{\R}\frac{\log w(x)}{1+x^2}\,dx
$$
converges. The trivial bound $\log w\le w$ shows that logarithmic integral of a Poisson-finite measure can diverge only to
$-\infty$. It will be convenient to call the set of all  measures  with finite logarithmic integral the Szeg\H{o} class $\sz$, i.e., 
$$
\sz = \Bigl\{\mu: \int_{\R}\frac{d\mu(x)}{1+x^2} + \int_{\R}\frac{|\log w(x)|}{1+x^2}\,dx < +\infty\Bigr\}.
$$
If $m\in \mathcal{N}(\C_+)$ and measure $\mu$ in \eqref{hg} is in Szeg\H{o} class, we can define
\begin{equation}\label{sd_h2}
\mathcal{K}_m = \log \Im m(i)-\frac{1}{\pi}\int_{\R}\frac{\log w(x)}{1+x^2}\,dx=\log\left(b+
\frac{1}{\pi}\int_{\R}\frac{d\mu}{1+x^2}\right)-\frac{1}{\pi}\int_{\R}\frac{\log w(x)}{1+x^2}\,dx\,.
\end{equation}\smallskip
One can use $b\ge 0$ and Jensen's inequality  to show that $\mathcal{K}_m\ge 0$. Notice that $\mathcal{K}_m=0$ if and only if $m$ is a constant with positive imaginary part. 

\smallskip

Let us introduce the class of Hamiltonians that characterizes  measures in Szeg\H{o} class.
If $\Hh$ is such that $\sqrt{\det\Hh} \notin L^1(\R_+)$, define
\begin{equation}\label{eq2}
\widetilde \K(\Hh) = \sum_{n=0}^\infty \left(\det \int_{\eta_n}^{\eta_{n+2}}\Hh(t)dt-4\right), \qquad \eta_n = \min\Bigl\{t \ge 0: \int_{0}^{t}\sqrt{\det{\Hh(s)}}\,ds = n\Bigr\}.
\end{equation}
Since the entries of $\Hh$
are locally integrable functions, the function $t \mapsto \sqrt{\det\Hh(t)}$ is also locally integrable on $\R_+$ and $\{\eta_n\}$ make sense. It is not difficult to check  (see Lemma \ref{sd_h1} in Appendix) that
$$
\det\int_{\eta_n}^{\eta_{n+2}} \Hh(t)\,dt \ge \left(\int_{\eta_n}^{\eta_{n+2}}\sqrt{\det\Hh(t)}\,dt\right)^2 = 4, \quad n \ge 0\,.
$$
This shows that the series in \eqref{eq2} contains only non-negative terms and hence its sum $\widetilde \K(\Hh) \in \R_+ \cup\{+\infty\}$ is well-defined but could be $+\infty$, in general. In Lemma \ref{sd_mg1}, we explain that $\widetilde\K(\Hh)$ can be rewritten in the form reminiscent of matrix $A_2$ Muckenhoupt condition. Roughly speaking, $\widetilde\K(\Hh)$ measures how fast the entries of $\Hh$ oscillate. In fact, we have $\widetilde\K(\Hh) = 0$ if and only if the Hamiltonian $\Hh$ is equivalent to a constant positive matrix, see Lemma \ref{sd_h1}. Notice that if the Hamiltonian is trivial then its determinant is zero and $\widetilde \K$ is undefined. Define the class $\hb$ of Hamiltonians by
$$
\hb = \Bigl\{\mbox{singular nontrivial }\Hh: \sqrt{\det \Hh} \notin L^1(\R_+), \; \widetilde \K(\Hh) < +\infty \Bigr\}.
$$
Here is the main result of the paper.
\begin{Thm}\label{t1} 
The spectral measure of a singular nontrivial Hamiltonian $\Hh$ on $\R_+$ belongs to the Szeg\H{o} class $\sz$ if and only if $\Hh \in \hb$. Moreover, we have 
\begin{equation}\label{ak_1}
c_1\K_m\le \widetilde\K(\Hh)\le c_2\K_me^{c_2\K_m},
\end{equation}
for some absolute positive constants $c_1$, $c_2$.
\end{Thm}
We emphasize that $\eqref{ak_1}$ is essentially sharp up to numerical values of $c_1$ and $c_2$. Indeed, for $\Hh$ such that $\widetilde\K(\Hh)   \le 1$, \eqref{ak_1} gives $\K_m\sim \widetilde\K(\Hh)$. Moreover, in Section 9 we present  two examples for both of which $\widetilde\K(\Hh)   > 1$. In the first example, we have $\K_m\sim \log (1+L)$ and $\widetilde\K(\Hh)\sim L$, where $L$ is arbitrarily large parameter. This shows that the exponent in the right hand side of \eqref{ak_1} can not be dropped. In the second example, we have $\K_m\sim L$ and $\widetilde\K(\Hh)\sim L$, where $L$ is again arbitrarily large parameter. Thus, the left bound in \eqref{ak_1} can not be improved. \medskip

The problem of controlling the entropy of the spectral measure for various differential operators has a long history and dates back at least to M.~Krein's  work \cite{Kr81} published in 1955. Quite recently, a large number of results that relate coefficients in differential or difference operators and spectral data were obtained (see, e.g., \cite{BBP}, \cite{B17}, \cite{Den02},\cite{Den06}, \cite{KS03}, \cite{KS09},  \cite{Mak_Polt}, \cite{NPVY}, \cite{Tep05},     and a book \cite{SimonDes}). Many of them can be considered as analogs of Szeg\H{o} theorem from the theory of polynomials orthogonal on the unit circle.  Our main theorem provides, perhaps, the most natural and far-reaching extension of this classical result. The following less general and a bit weaker version of Theorem \ref{t1} has been proved in \cite{BD2017}.
\begin{Thm}{\rm (Bessonov-Denisov, \cite{BD2017})}\label{bd-2017}
An even measure $\mu$ belongs to the Szeg\H{o} class $\sz$ if and only if some (and then every)  Hamiltonian 
$\Hh = \left(\begin{smallmatrix} h_1 &0\\0 &h_2 \end{smallmatrix}\right)$ 
generated by $\mu$ is such that $\sqrt{\det \Hh} \notin L^1(\R_+)$ and 
\begin{equation}\label{eq2_o2}
\widetilde \K(\Hh) = \sum_{n = 0}^{+\infty}  \left(\int_{\eta_n}^{\eta_{n+2}}h_1(s)\,ds \cdot \int_{\eta_n}^{\eta_{n+2}}h_2(s)\,ds - 4\right) < \infty,
\end{equation}
where $\{\eta_n\}$ are given by \eqref{eq2}. 
Moreover, we have $\widetilde\K(\Hh) \le c\K_me^{c\K_m}$ and $\K_m \le c\widetilde\K(\Hh) e^{c\widetilde\K(\Hh)}$ for an absolute constant $c$. 
\end{Thm}

A characterization of Krein strings for which the spectral measure has finite logarithmic integral has been given in \cite{BD2017} as well. That was an immediate corollary of Theorem~\ref{bd-2017}. Other spectral theoretic applications of Theorems \ref{t1} and \ref{bd-2017} can be found in \cite{B2018},\cite{B18a},\cite{EK19},\cite{EKK19},\cite{NK19}.

\medskip

Some proofs in \cite{BD2017} relied on the fact that diagonal matrices (arising from diagonal Hamiltonians) commute, which forces us to find a 
different argument for the proof of Theorem~\ref{t1} in full generality. We also want to emphasize here that the method used in our proof does not involve any ``sum rules'', which often times is the basis for other proofs found in the literature. We outline the main steps of the proof in Section \ref{s3}.

\medskip

The Szeg\H{o} class proved to be important in mathematical physics, in particular, in the scattering theory of wave propagation. For example, in \cite{Den02b}, strong wave operators for a one-dimensional Dirac system with a $L^2(\R_+)$-potential were expressed in terms of the Szeg\H{o} function of the spectral measure.  The main result of  \cite{B2018} shows that regularized version of strong wave operator for a one-dimensional Dirac system exists and is complete under the single assumption that the spectral measure belongs to the Szeg\H{o} class $\sz$. Using Theorem \ref{t1}, we describe such Dirac systems below.
\begin{Cor}\label{t2} 
Let $\mu$ be the spectral measure of the Dirac operator $\mathcal{D}_{V}$ on $\R_+$, defined by
\begin{equation}\label{sd_dirac}
\mathcal{D}_{V}:\left(\begin{smallmatrix}f_1\\f_2 \end{smallmatrix}\right) \mapsto J \frac{d}{dt} \left(\begin{smallmatrix}f_1\\f_2 \end{smallmatrix}\right) + V(t)\left(\begin{smallmatrix}f_1\\f_2 \end{smallmatrix}\right),  \quad t \in \R_+,\quad f_2(0)=0,
\end{equation}
with a real-valued locally summable $2 \times 2$ potential $V = V^*$  which satisfies condition $\, \trace V = 0$. Then, $\mu \in \sz$ if and only if $N_0^*N_0 \in  \hb$, where $N_0$ solves $J N_0'(t) + V(t)N_0(t) = 0$, $N_0(0) = I_{2\times 2}$, $t \in \R_+$.
\end{Cor}
One version of classical Krein-Wiener completeness theorem says that the future subspace of a  Gaussian stationary process is not determined by its past subspace if and only if the spectral measure of the process belongs to the Szeg\H{o} class, see, e.g., \cite{ibr}. Very interesting direction for further research is to find probabilistic applications of Theorem \ref{t1}. We mention two  papers \cite{Alpay}, \cite{Dzh} related to the subject. 

\medskip

A few months after the current manuscript was posted on arXiv, the authors received a note from Peter Yuditskii in which the logarithmic integral of a quantity closely connected to spectral measure was expressed via the integral of elements of Hamiltonian, written in a special form. It is of interest to relate this ``sum rule'' to estimates obtained in this work. 

\medskip

Here is an outline of the paper. In the second section, we give more detail about canonical systems. In the third section, we explain the main steps of the proof of the main result, Theorem~\ref{t1}. Section \ref{s41} contains some examples relevant  to Theorem \ref{t1}. It is followed by sections which contain different parts of the proof. In the Appendix, we collect auxiliary results used in the main text.

\medskip

\subsection{Notation}
\begin{itemize}
\item ${\rm SL}(2,\R)$ denotes the set of real $2\times 2$ matrices with unit determinant.
\item If $A$ is $d\times d$ matrix, $\|A\|$ stands for operator norm in $\C^d$.
\item For $p\ge 1$, let us denote by $L^p$ the set of $2\times 2$ matrix-valued functions $V$ on $\R_+$ such that 
$$
\|V\|^p_{L^p}\dd \int_{\mathbb{R}_+}\|V(t)\|^pdt<\infty.
$$
Let $L^1+L^2$ denote the set of sums $V=V_1+V_2$ equipped with the norm 
$$
\|V\|_{1,2}\dd \|V\|_{L^1+L^2}\dd \inf \{\|V_1\|_{L^1}+\|V_2\|_{L^2}: V=V_1+V_2\}.
$$ 
Similar notation will be used for scalar functions.
\item The symbol $C$ denotes the absolute constant which can change the value from formula to formula. If we write, e.g., $C(\alpha)$, this defines a positive function of parameter $\alpha$.
\item For two non-negative functions $f_1$ and $f_2$, we write $f_1\lesssim f_2$ if  there is an absolute
constant $C$ such that $f_1\le Cf_2$ for all values of the arguments of $f_1$ and $f_2$. We define $\gtrsim$
similarly and say that $f_1\sim f_2$ if $f_1\lesssim f_2$ and
$f_2\lesssim f_1$ simultaneously. If $|f_3|\lesssim f_4$, we will write $f_3=O(f_4)$. 
\item Given any interval $I\subset \R$ and $f\in L^1(I)$, we define
$$
\langle f\rangle_I= \frac{1}{|I|}\int_I fdx\,.
$$
\item Entire function $E$ is called Hermite-Biehler function if it has no zeroes in $\C_+$ and 
$$
|E^\sharp(z)|\le |E(z)|, \quad z\in \C_+,
$$
with 
$$
E^\sharp(z)\dd \overline{E(\bar{z})}.
$$
\item If $S$ is a set, $\chi_S$ denotes the characteristic function of $S$.
\item We sometimes use symbol $\K_{\Hh}(0)$ instead of $\K_m$. The reader should be aware that  these two quantities are identical by definition. Notation $\K_{\Hh}(0)$ will be explained in the next section.
\item $\mathbb{Z}_+=\{0,1,\ldots\}.$
\item If $A$ is self-adjoin matrix, we denote  its smallest and largest eigenvalues by $\lambda_{\min}(A)$ and $\lambda_{\max}(A)$, respectively.
\item For $a>0$, we define
\begin{equation}
\log^+ a=
\begin{cases}\label{sd_vi1}
\log a, & a\ge 1,\\
0,&a\in (0,1),
\end{cases}
\qquad 
\log^- a=\begin{cases}
-\log a, & a\in (0,1],\\
0,&a>1.
\end{cases} 
\end{equation}
so $\log^+a\ge 0$, $\log^-a\ge 0$, and $\log a=\log^+a-\log^-a$.
\end{itemize}

\medskip

\section{Preliminaries on canonical Hamiltonian systems}\label{s2}
In this section, we collect some definitions and known results that will be used later in the text. In fact,
we almost literally repeat the content of Section 1 in \cite{BD2017} and Section 2 in \cite{B2018}. See monographs \cite{dbbook}, \cite{Remlingb},
\cite{Romanov} for the classical theory of de Branges systems.

\medskip

\subsection{Two results on canonical systems}

Later in the text, we will need two classical results from the spectral theory of canonical Hamiltonian systems. Given a
Hamiltonian $\Hh$ on $\R_+$, define the function
$$
\xi_\Hh: t \mapsto \int_{0}^{t}\sqrt{\det \Hh(s)}\,ds,
\qquad t \in \R_+.
$$
Since $2\sqrt{\det\Hh(s)} \le \trace \Hh(s)$ for all $s \ge 0$, the function $s \mapsto \sqrt{\det\Hh(s)}$ is integrable on compact
subsets of~$\R_+$. In particular, $\xi_\Hh$ is correctly defined and absolutely continuous function on $\R_+$. In
the case when $\sqrt{\det\Hh} \notin L^1(\R_+)$, one can define the function
$$
\eta_{\Hh}: t \mapsto \min\Bigl\{r \ge 0: \;\; t = \int_{0}^{r}\sqrt{\det\Hh(s)}\,ds\Bigr\}, \qquad t \in \R_+.
$$
Observe that $\eta_{\Hh}(n) = \eta_n$ for $\eta_n$ in \eqref{eq2}. For $s > 0$, denote by $\E_s$ the linear space of functions with
smooth Fourier transform supported on $(0, s)$. The following theorem is a consequence of results by
M.~Riesz, S.~Mergelian, and M.~Krein, see Proposition 2.5 in \cite{B2018}.
\begin{Thm}\label{p1}
Let $\Hh$ be a singular nontrivial Hamiltonian on $\R_+$, and let $\mu$ be its spectral
measure. If $\E_s$ is not dense in $L^2(\mu)$ for some $s > 0$, then $\xi_\Hh(t) \ge s$ for some $t \ge 0$.
\end{Thm}
Next result is usually referred to as the Krein-Wiener completeness theorem. See Section~4.2
in \cite{DymMcKean} or Theorem A.6 in \cite{Den06} for the proof.
\begin{Thm}\label{p2} 
Let $\mu$ be a Poisson-finite measure on $\R$. Then $\mu \in \sz$ if and
only if $\bigcup_{s>0}\E_s$ is not dense in $L^2(\mu)$.
\end{Thm}

\noindent {\bf Remark.} Our main result, Theorem \ref{t1}, complements Krein-Wiener's theorem by giving yet another criterion for completeness.

\medskip

\subsection{Bernstein-Szeg\H{o} approximation, entropy function of a Hamiltonian,  ${\rm SL}(2,\R)$ invariance} 
Let $\Hh$ be a singular nontrivial Hamiltonian on $\R_+$. For every $r > 0$, define $\Hh_r$ to be the Hamiltonian 
$t \mapsto \Hh(t + r)$ defined on $\R_+$. Let $m_r$, $\mu_r$, $b_r$, $a_r$ denote the Titchmarsh-Weyl function of $\Hh_r$, its spectral measure, and the coefficients in the Herglotz representation \eqref{hg}
for $m_r$. Define
\begin{align}
\I_{\Hh}(r) &= \Im m_r(i) = \frac{1}{\pi}\int_{\R}\frac{d\mu_r(x)}{1+x^2}+b_r,\notag\\
\Rr_{\Hh}(r) &= \Re m_r(i) = a_r, \label{sd_h66}\\
\J_{\Hh}(r) &= \frac{1}{\pi}\int_{\R}\frac{\log w_r(x)}{1+x^2}\,dx,\notag
\end{align}
where $\mu_r = w_r\,dx + \mu_{r, \mathbf{s}}$ is the decomposition of $\mu_r$ into the absolutely continuous and singular parts. In the case when $\mu_r \notin \sz$ for some $r \ge 0$, we set $\J_\Hh(r) = -\infty$. The
entropy function of $\Hh$ is introduced as follows
\begin{equation}\label{sd_h17}
\K_\Hh(r) = \log\I_\Hh(r) - \J_\Hh(r), \qquad r \ge 0.
\end{equation}
Notice that $\K_\Hh(0)=\K_m$, where $\K_m$ was defined in \eqref{sd_h2}. Since $b_r \ge 0$, Jensen's inequality implies $\K_\Hh(r) \ge 0$. Next, consider the Hamiltonian
\begin{equation}\label{cook59}
\widehat\Hh_r(t) = 
\begin{cases}
\Hh(t),  &t \in [0,r),\\
\Bigl(\begin{smallmatrix} c_1(r) & c(r)\\ c(r) & c_2(r)\end{smallmatrix}\Bigr), &t \in [r,+\infty),
\end{cases}
\end{equation}
where $c_1(r) = 1/\I_\Hh(r)$, $c(r) = \Rr_{\Hh}(r)/\I_{\Hh}(r)$, $c_2(r) = (\I^2_\Hh(r) + \Rr^2_{\Hh}(r))/\I_{\Hh}(r)$. The Hamiltonian $\widehat\Hh_r$
coincides with $\Hh$ on $[0, r)$ and is constant on $[r, +\infty)$.  We call $\widehat\Hh_r$ the Bernstein-Szeg\H{o} approximation to~$\Hh$. Some properties of the functions $\K_{\Hh}$, $\I_{\Hh}$, $\Rr_\Hh$ are collected in the following two lemmas.
\begin{Lem}\label{l1} 
Let $\Hh = \sth$ be a singular nontrivial Hamiltonian on $\R_+$ and let $\mu = w dx + \mu_{\mathbf{s}}$
be the spectral measure of $\Hh$. Assume that $\mu \in \sz$. Then, for every $r \ge 0$ the measure $\mu_r =
w_r \, dx + \mu_{r,\mathbf{s}}$ belongs to $\sz$ and
\begin{itemize}
\item[$(a)$] $\K_{\Hh}(0) = \K_{\widehat\Hh_r}(0) + \K_{\Hh}(r)$,
\item[$(b)$] $\lim_{r \to +\infty}\K_{\Hh}(r) = 0$, $\lim_{r \to +\infty}\K_{\widehat\Hh_r}(0) = \K_{\Hh}(0)$.
\end{itemize}
If, moreover, $\det \Hh = 1$ a.e.  on $\R_+$, then $\K_{\Hh},\I_{\Hh}$, and $\Rr_{\Hh}$ are absolutely continuous on $\R_+$ and
\begin{itemize}
\item[$(c)$] $   \displaystyle \K'_{\Hh} = \left(2 - \I_\Hh h_1 - \frac{1}{\I_\Hh h_1} \right)- \frac{(\Rr'_\Hh/\I_\Hh)^2}{4\I_\Hh h_1}$,
\item[$(d)$] $   \displaystyle \I'_\Hh/\I_\Hh = \left(\I_\Hh h_1 - \frac{1}{\I_\Hh h_1 } \right)- \frac{(\Rr'_\Hh/\I_\Hh)^2}{4\I_\Hh h_1}$,
\item[$(e)$] $\Rr'_\Hh/\I_\Hh = 2\Rr_\Hh h_1 - 2h$,
\end{itemize}
almost everywhere on $\R_+$.
\end{Lem}
\beginpf For items  $(a)$ and $(b)$, see Lemma 2.3 in \cite{B2018} and its proof therein (Appendix I in \cite{B2018}). Identities
$(c)$-$(e)$ are equivalent to formulas $(39)$-$(41)$ in \cite{B2018} after elementary algebraic manipulations are  performed. \qed

\medskip

\noindent {\bf Remark.} In the case of diagonal $\Hh$, identities $(a)$-$(e)$ can be found in Lemma~2.5 and Lemma~2.7 in \cite{BD2017}.\bigskip

\noindent Recall that SL$(2,\R)$ is related to fractional linear transformations that leave $\C_+$ invariant, i.e.,
$$
A=
\begin{pmatrix}
a &b\\c & d
\end{pmatrix} 
\quad \Longleftrightarrow \quad {\rm Mob}_A(z)\dd \frac{az+b}{cz+d}, \quad a,b,c,d\in \mathbb{R}, \quad ad-cb=1.
$$
Given any $A\in {\rm SL}(2,\R)$, we define conjugation of $\Hh$ by $A$ as follows: 
$$
\Hh_A = A^*\Hh A\,.
$$
The spectral measure of $\Hh_A$ will be denoted by $\mu_A$.
Next lemma proves that both $\widetilde \K({\Hh})$ and $\K_{\Hh}$ are invariant under conjugation and under linear fractional transform, respectively.
\begin{Lem}\label{li2} 
Let $\Hh$ be a singular nontrivial Hamiltonian on $\R_+$, and let $\mu$ be its spectral measure.
 Then, 
 \begin{itemize}
\item[(a)] $\Hh\in \hb$ if and only if $\Hh_A \in \hb,$ 
\item[(b)] $\mu \in \sz$ if and only if $\mu_A \in \sz$. 
 \end{itemize}
 Moreover,
$\widetilde\K(\Hh) = \widetilde \K(\Hh_A)$ and $\K_{\Hh}(0) = \K_{\Hh_A}(0)$ whenever these quantities are finite.
\end{Lem}

\beginpf Since $\det \Hh_A(t) = \det \Hh(t)$ for a.e.  $t \in \R_+$ and $A$ is $t$-independent, we see that $\Hh\in\hb$ if and only if $\Hh_A\in \hb$, and, moreover, $\widetilde\K(\Hh) = \widetilde \K(\Hh_A)$. This proves $(a)$.

\smallskip

To show $(b)$, we first compute Titchmarsh-Weyl function for $\Hh_A$. Suppose $M$ and $M_A$ are the solutions of \eqref{eq1} for the Hamiltonians
$\Hh$ and $\Hh_A$, respectively. We are claiming that 
\begin{equation}\label{amd_1}
M_A = A^{-1} M A\,.
\end{equation}
Indeed, by definition of $M$ we have $JM'=z\Hh M$, $M(0,z)=I_{2\times 2}$. This implies
$$
(A^*JA)(A^{-1}MA)'=z(A^*\Hh A)(A^{-1}MA),\quad A^{-1}M(0,z)A=I_{2\times 2}.
$$
To prove \eqref{amd_1}, we only need to notice that $A^*JA=J$ by Lemma \ref{sd_h7}. For $A=\left(\begin{smallmatrix}a&b\\c&d\end{smallmatrix}\right)$, we have
\begin{equation}
M_A =
\begin{pmatrix}
d & -b\\
-c & a
\end{pmatrix}
\begin{pmatrix}
\Theta^+ & \Phi^+\\
\Theta^- & \Phi^-
\end{pmatrix}
\begin{pmatrix}
a & b\\
c & d
\end{pmatrix}.\label{sd_h62}
\end{equation}
Taking $\omega = 0$ in \eqref{mf} for $\Hh_A$, we get
\begin{equation}\label{li16}
m_A(z) 
=\lim_{t\to +\infty}
\frac{
(-c\Theta^+(z,t) + a\Theta^-(z,t))b + (-c\Phi^+(z,t) + a\Phi^-(z,t))d}{ (-c\Theta^+(z,t) + a\Theta^-(z,t))a + (-c\Phi^+(z,t) + a\Phi^-(z,t))c}
=\frac{d m(z) + b}{c m(z) + a}, \qquad z \in \C_+.\notag
\end{equation}
It remains to note that $\Im m_A=\frac{\Im m}{|cm+a|^2}$, hence
\begin{align*}
\K_{\Hh_A}(0) 
&= \log \frac{\Im m(i)}{|c m(i) + a|^2} - \frac{1}{\pi} \int_{\R}\frac{\log w(x) -\log (|cm(x) + a|^2)}{x^2+1}\,dx,\\
&= \log \Im m(i) - \frac{1}{\pi} \int_{\R}\frac{\log w(x)}{x^2+1}\,dx = \K_{\Hh}(0),
\end{align*}
where we used Lemma \ref{sd_h13} in Appendix. \qed 

\medskip

\noindent The Hamiltonian dual to $\Hh$ is defined by conjugating with $A=J$, i.e.,
$$
\Hh_d\dd J^*\Hh J.
$$
Notice that $(\Hh_d)_d=\Hh$. Lemma \ref{li2} yields the following corollary (see also Lemma 3 in \cite{B2018}).   
\begin{Cor}\label{c1}
Let $\Hh$ be a singular nontrivial Hamiltonian, $\mu_d$ denote the spectral measure of $\Hh_d$, and $m_d$ denote the Titchmarsh-Weyl function of $\Hh_d$. Then, $\mu \in \sz$ if and only if  $\mu_d\in \sz$. Moreover, we have $\K_\Hh = \K_{\Hh_d}$ and $m_d = -1/m$.
\end{Cor}
\beginpf The first part of the statement follows from Lemma \ref{li2} by taking $A = J$. Formula \eqref{li16} shows that $m_d = -1/m$. \qed

\medskip

\noindent {\bf Remark.} In \eqref{sd_h2}, the definition  of $\K_m$, we evaluate $\Im m$ at  $z_0=i$ and the Poisson kernel inside the integral is evaluated at the same point. Changing this reference point results in the whole family of entropies indexed by parameter $z_0\in \C_+$. Clearly, if entropy at one point is finite, it is finite at any other point. In this paper, we do not study how our main result can be modified (i.e., how constants in two-sided estimates in Theorem \ref{t1} depend on $z_0$) but we believe this is a promising direction. \bigskip

\section{Main steps in the proof of Theorem \ref{t1}}\label{s3}

In this short section, we explain the structure of the proof of Theorem \ref{t1}. The following special factorization of Hamiltonian lies at the core of our approach.

\medskip

\noindent {\bf Definition.} Suppose $\mathfrak{q},\mathfrak{v}_1,\mathfrak{v}_2$ are three non-negative parameters. Let $\Hh$ be a Hamiltonian which satisfies $\det \Hh=1$ for a.e. $t\ge 0$. We will say that $\Hh$ admits $(\mathfrak{q},\mathfrak{v}_1,\mathfrak{v}_2)$ -- factorization if $\Hh=G^*QG$ for some $2\times 2$ matrix-valued functions $Q$, $G$ with real entries such that
 \begin{equation}\label{sd_list}\left\{
\begin{array}{cl}
(a) & Q\ge 0,\, \det Q=1\,\, {\rm a.e.\,  on} \,\,\mathbb{R}_+,\\
(b) & \|\trace  Q-2\|_{L^1(\mathbb{R}_+)}\le \mathfrak{q},\\
(c) & G \,\,{\rm  is \,\, absolutely \,\,  continuous\,\, and } \,\, \det G=1\,\, {\rm  on} \,\, \mathbb{R}_+,\\
(d) & G'=JVG  \,\, {\rm  for   \,\, some}  \,\, V=V^*, V=V_1+V_2  \,\, {\rm  with} \\
\,\, &V_1\in L^1, V_2\in L^2 \,\,{\rm such \,\,  that} \,\,\|V_1\|_{L^1} \le \mathfrak{v}_1 \,\,{\rm and }   \,\,\|V_2\|_{L^2} \le \mathfrak{v}_2.
\end{array}\right.
\end{equation}
\noindent {\bf Remark.} This $(\mathfrak{q},\mathfrak{v}_1,\mathfrak{v}_2)$--factorization is not unique, in general. Since $Q\ge 0$ and $\det Q=1$, we get $\trace Q\ge 2$. Note that the matrix-valued function $V$ in this definition necessarily has real entries. Equation $G'=JVG$ gives
$$
\det G(t)=\det G(0)\cdot \exp\left(\int_0^t \trace (JV(\tau))d\tau\right).
$$
For $2\times 2$ real matrices, the condition $\trace(JV)=0$ is equivalent to $V$ being symmetric. Therefore, $V$ being symmetric and $\det G(0)=1$ already imply $\det G(t)=1$ for all $t$. The factor $G$ can be regarded as ``slow''  and factor $Q$ can be regarded as the ``fast'' one. Indeed, elements of $G$ are absolutely continuous and elements of $Q$ are only locally integrable. On the other hand, $G(t)$ can grow infinitely when $t\to+\infty$ although $Q(t)$ is ``close to $I_{2\times 2}$'' at infinity as follows from $(a)$ and $(b)$. \smallskip

\noindent {\bf Remark.} In the definition of Hamiltonians that admit factorization, we didn't specify $G(0)$. This was done intentionally. In fact, given $G,Q$ and $\Hh=G^*QG$, we can take $\Hh_{G^{-1}(0)}=(GG^{-1}(0))^*Q(GG^{-1}(0))$. The parameters $Q,V$ in the factorization of $\Hh_{G^{-1}(0)}$ can be chosen the same as in that of $\Hh$  and $\K_\Hh(0)=\K_{\Hh_{G^{-1}(0)}}(0)$, $\widetilde \K(\Hh)=\widetilde \K(\Hh_{G^{-1}(0)})$ by Lemma \ref{li2} (and we will later prove that these quantities are in fact finite).

\medskip

\noindent {\bf Definition.} The class $\fc$ (shorthand for ``finally constant'') is the set of singular nontrivial Hamiltonians $\Hh$ on $\R_+$ such that $\Hh = A$ on $[\ell, +\infty)$ for some $\ell \ge 0$ and a constant positive matrix $A$. Parameter $\ell$ and matrix $A$ might depend on $\Hh$.
 
\medskip

Theorem \ref{t1} follows from five theorems formulated below. Their proofs are given in Sections~\ref{s5}-\ref{s9}.
\begin{Thm}\label{l8}
Assume that 
\begin{equation}
\label{eq15_1}
c_1\K_m\le \widetilde\K(\Hh)\le c_2\K_me^{c_2\K_m}
\end{equation}
holds for every Hamiltonian $\Hh\in \fc$ such that $\det \Hh = 1$ almost everywhere on $\R_+$. Then, the conclusions of Theorem \ref{t1} follow, i.e.,
\begin{itemize}
\item[$(a)$] The spectral measure $\mu\in \sz$ of a singular nontrivial Hamiltonian $\Hh$ belongs to $\sz$ if and only if $\Hh\in \hb$.
\item[$(b)$] Moreover, \eqref{eq15_1} holds for all $\Hh \in \hb$ with the same constants $c_1,c_2$, and $c_3$. 
\end{itemize}
\end{Thm}
\begin{Thm}\label{sd_h32}
Let $\Hh$ be a singular nontrivial Hamiltonian which satisfies $\det \Hh=1$ almost everywhere on $\R_+$, and let $\mu$ be its spectral measure. If $\mu\in \sz$, then $\Hh$ admits $(\mathfrak{q}, \mathfrak{v}_1, \mathfrak{v}_2)$--factorization with $\mathfrak q\lesssim \K_{\Hh}(0), \mathfrak{v}_1\lesssim \K_\Hh(0),$ and $\mathfrak{v}_2\lesssim \sqrt{\K_\Hh(0)}$.
\end{Thm}

\begin{Thm}\label{t42}
Let $\Hh$ be a singular nontrivial Hamiltonian on $\R_+$ and let $\mu$ be its spectral measure. If $\Hh$ admits  $(\mathfrak{q}, \mathfrak{v}_1, \mathfrak{v}_2)$ -- factorization, then $\mu\in \sz $. Moreover,  
$\K_\Hh (0)\lesssim \min\{\mathfrak{v}_1,\mathfrak{v}_1^2\}+\mathfrak{v}_2^2+\mathfrak{q}$.
\end{Thm}
\begin{Thm}\label{t5}
Suppose that $\Hh$ is a Hamiltonian on $\R_+$ allowing $(\mathfrak{q}, \mathfrak{v}_1, \mathfrak{v}_2)$--factorization. Then $\Hh \in \hb$ and we have $\widetilde \K(\Hh)  \le c(\mathfrak{q} + \mathfrak{q}^2 + \mathfrak{v}_1^2 + \mathfrak{v}_2^2)e^{c\mathfrak{v}_1 + c\mathfrak{v}_2}$ for an absolute constant $c$. 
\end{Thm}   
\begin{Thm}\label{t6}
Suppose that $\Hh \in \hb$ and $\det \Hh = 1$ for almost all $t\in\R_+$. Then $\Hh$ admits $(\mathfrak{q}, \mathfrak{v}_1, \mathfrak{v}_2)$-factorization. Moreover, we have $\mathfrak{q} \lesssim  \widetilde \K({\Hh}) $, $\mathfrak{v}_1 \lesssim \widetilde \K({\Hh})$, and $\mathfrak{v}_{2}^{2} \lesssim \widetilde \K({\Hh})$. 
\end{Thm}

\medskip

Assuming Theorems \ref{l8}--\ref{t6} are proved, we can easily finish the proof of the main result.

\medskip

\noindent{\bf Proof of Theorem \ref{t1}.} By Theorem \ref{l8}, it suffices to show that 
\begin{equation}\label{eq15bis}
c_1\K_m\le \widetilde\K(\Hh)\le c_2\K_me^{c_2\K_m},
\end{equation}
for every Hamiltonian $\Hh$ with unit determinant, which belongs to class ${\rm FC}$. Take such $\Hh$. Combining
Theorem \ref{t42} and Theorem \ref{t6}, we see that $\Hh$ admits $(\mathfrak{q},\mathfrak{v}_1,\mathfrak{v}_2)$-factorization. Moreover, we
get the estimates
$$
\K_\Hh (0)\lesssim  \min\{\mathfrak{v}_1,\mathfrak{v}_1^2\}+ \mathfrak{v}_2^2+  \mathfrak{q}, \qquad \mathfrak{q} \lesssim \widetilde \K({\Hh}),\qquad \mathfrak{v}_{1} \lesssim  \widetilde \K({\Hh}),\qquad \mathfrak{v}_2^2 \lesssim  \widetilde \K({\Hh}),
$$
so $\K_{\Hh}(0)\lesssim \widetilde \K(\Hh)$. From Theorem \ref{t5} and Theorem \ref{sd_h32}, we get
$$
 \widetilde \K(\Hh) \le c(\mathfrak{q} + \mathfrak{q}^2 + \mathfrak{v}_1^2 + \mathfrak{v}_2^2)e^{c\mathfrak{v}_1 + c\mathfrak{v}_2}, \qquad
 \mathfrak q\lesssim \K_{\Hh}(0),\qquad \mathfrak{v}_1\lesssim \K_\Hh(0), \qquad \mathfrak{v}_2\lesssim \sqrt{\K_\Hh(0)},
$$
so $\widetilde\K({\Hh})\le c\K_{\Hh}(0)e^{c\K_{\Hh}(0)}$ and we have \eqref{eq15bis}. \qed

\medskip

\noindent {\bf Proof of Corollary \ref{t2}}. It is known that the spectral measure of Dirac system \eqref{sd_dirac} coincides with the spectral measure of the canonical Hamiltonian system generated by the Hamiltonian $\Hh = N_0^*N_0$, see details in \cite{B2018}. Thus, the application of Theorem \ref{t1} gives the corollary. \qed

\bigskip

\section{Hamiltonians in class $\hb$, matrix-valued $A_2$-condition, and some examples}\label{s41}
The diagonal Hamiltonians in class $\hb$ have been thoroughly studied in \cite{BD2017}. If we assume that $\Hh=\left(\begin{smallmatrix}h_1&0\\0&h_1^{-1}\end{smallmatrix}\right)$, i.e., if $\Hh$ is diagonal and $\det \Hh=1$ for a.e. $t\in \R_+$, then condition $\Hh\in \hb$ reads as follows:
$$
\widetilde\K(\Hh)=\sum_{n=0}^\infty \left(\left(\int_n^{n+2} h_1dt\right)\left(\int_n^{n+2} h_1^{-1}dt\right)-4\right)<\infty.
$$
In \cite{BD2017}, the class of functions $h_1$ that satisfy this condition was denoted by $A_2(\R_+,\ell^1)$ in analogy to the standard Muckenhoupt $A_2$ condition on the weights. We recall that, given non-negative matrix-valued function $W$ defined on $\R$, the matrix $A_2$ Muckenhoupt characteristics of $W$ is defined by (see, e.g., \cite{treil_volberg})
$$
[W]_{A_2}\dd \sup_{I}\left\|   \langle W\rangle_I^{1/2}\langle W^{-1}\rangle_I^{1/2}\right\|,
$$
where the supremum is taken over all intervals $I$ in $\R$. To see the connection with our condition \eqref{eq2}, we need the following lemma. 

\begin{Lem}\label{sd_mg1}
Suppose $H$ is $2\times 2$ nonnegative matrix-valued function defined on $I\dd [a,b]$ and $H$  satisfies $\det H=1$ for a.e. $t\in I$. Then, 
$$
\|\langle H\rangle_I^{1/2}\langle H^{-1}\rangle_I^{1/2}\|={\det}^{1/2} \langle H\rangle_I.
$$
In particular, for every $\Hh\in \hb$ that satisfies $\det \Hh=1$, we have
\begin{equation}\label{eq91}
\widetilde\K(\Hh)=4\sum_{n=0}^\infty \Bigl( \|\langle \Hh\rangle_{[n,n+2]}^{1/2}\langle \Hh^{-1}\rangle_{[n,n+2]}^{1/2}\|^2-1\Bigr) < +\infty.
\end{equation}
\end{Lem}
\beginpf
By a change of variables, we can assume that $I=[0,1]$.  Let $\Omega= \langle H\rangle_I^{1/2}$. Since $\det H=1$ and $H=H^*$, we have 
$$
\langle H^{-1}\rangle_I=\int_0^1 H^{-1}dt\stackrel{\eqref{sd_h7}}{=}(-J)\left(\int_0^1 Hdt\right) J=(-J)\langle H\rangle_I J.
$$
Notice that 
$$
\left((-J)\left(\int_0^1 Hdt\right) J\right)^{1/2}=(-J)\Omega J,
$$
as can be checked directly. Then, since $\|A\|=\|A^*\|$ for every matrix $A$, we can write
$$
\|\Omega(-J)\Omega J\|=\|J\Omega J\Omega \|=\|\Omega J\Omega \|,
$$
because $J$, being the unitary matrix, preserves the norm.
Notice that for all positive $2\times 2$ matrices $\Omega=\left(\begin{smallmatrix} a_1&a\\a&a_2\end{smallmatrix}\right)$, we have an identity
$$
\|\Omega J\Omega\|=\det \Omega,
$$
which follows from the formula
$$
\Omega (iJ) \Omega=\left( \begin{smallmatrix} 0 & i(a^2-a_1a_2)\\-i(a^2-a_1a_2)&0\end{smallmatrix}   \right)
$$
and an observation that the last self-adjoint matrix has eigenvalues $\pm (a_1a_2-a^2)=\pm \det \Omega$. \qed

\medskip

We will call the class of weights satisfying \eqref{eq91} the matrix-valued $A_2(\R_+,\ell^1)$ class. The following lemma asserts that 
the diagonal elements of mappings in the matrix-valued class $A_2(\R_+,\ell^1)$ belong to the scalar class $A_2(\R_+,\ell^1)$.

\begin{Lem}\label{l42}
Let  $\Hh = \sth$ belong to $\hb$ and $\det \Hh=1$ a.e. on $\R_+$. Then, we have
$$
\sum_{n=0}^\infty\left(\left(\int_{n}^{n+2}h_{1}dx\right)\left(\int_{n}^{n+2}h_{1}^{-1}dx\right)-4\right)\leq \widetilde{\mathcal{K}}(\Hh).
$$
Similar bound holds for $h_{2}$.
\end{Lem}
\beginpf
For every interval $I$, we have by Cauchy-Schwarz inequality
$$
\langle h_1 \rangle_{I} \langle h_1^{-1} \rangle_{I}+\langle h \rangle_{I}^2\le\langle h_1 \rangle_{I} \langle (h^2+1)h_1^{-1} \rangle_{I}.
$$
Recall that $h_1h_2-h^2=1$  so $h_2=(1+h^2)h_1^{-1}$. Then, we can rewrite the last bound as
$$
\langle h_{1} \rangle_{I} \langle h_{1}^{-1} \rangle_{I}+\langle h \rangle_{I}^2\le \langle h_1 \rangle_{I} \langle h_2 \rangle_{I}, \quad \langle h_1 \rangle_{I} \langle h_1^{-1} \rangle_{I}\le \langle h_1 \rangle_{I} \langle h_2 \rangle_{I}-\langle h \rangle_{I}^2.
$$
Taking $I=[n,n+2]$, subtracting $1$ from both sides and summing in $n$ finishes the proof. \qed

\medskip

In the case of diagonal Hamiltonians, the proofs of Theorems \ref{t5}, \ref{t6} are much easier because they can be reduced to considerations of scalar functions. For instance, the following lemma solves the problem of existence of $(\mathfrak{q},\mathfrak{v}_1,\mathfrak{v}_2)$ -- factorization for diagonal Hamiltonians.
\begin{Lem}\label{l43}
A function $h$ on $\R_+$ belongs to $A_2(\R_+, \ell^1)$ if and only if there exist functions
$q$, $v$ on $\R_+$ such that 
\begin{itemize}
\item[$(a)$] $q > 0$ almost every where on $\R_+$, $q+q^{-1}-2\in L^1(\R_+)$,   
\item[$(b)$] $v$ is real-valued, $v \in L^1(\R_+)+L^2(\R_+)$,
\item[$(c)$] $h(t) = q(t)\exp\left(\int_0^t v(\tau)\, d\tau\right)$, $t \in \R_+$. 
\end{itemize}
\end{Lem}

\noindent While this lemma could be proved by means of elementary function theory, its proof is rather complicated. For completeness, we give a short proof based on Theorems \ref{t5}, \ref{t6}. 

\medskip

\beginpf Suppose that $h \in A_2(\R_+, \ell^1)$ and consider the Hamiltonian $\Hh = \left(\begin{smallmatrix}h & 0 \\ 0 & 1/h\end{smallmatrix}\right)$. We have $\Hh \in \hb$. An inspection of the proof of Theorem \ref{t6} shows that $\Hh$ admits $(\mathfrak{q},\mathfrak{v}_1,\mathfrak{v}_2)$ -- factorization of the form $\Hh=G^*QG$ with parameters $G$, $Q$ such that 
$$
Q=\begin{pmatrix}q & 0\\ 0 &q^{-1}\end{pmatrix}, \qquad G' = JVG, \quad G(0) = I_{2 \times 2}, \qquad 
V=
\begin{pmatrix}
0 & -v/2\\ -v/2 & 0
\end{pmatrix}, 
$$
where $q$ satisfies $(a)$ and $v$ satisfies $(b)$. Solving equation $G' = JVG$, $G(0) = I_{2 \times 2}$, we get  
$$
G=\begin{pmatrix}
e^{\frac{1}{2}\int_0^t vd\tau} & 0\\
0 &  e^{-\frac{1}{2}\int_0^t vd\tau}
\end{pmatrix}, 
\qquad 
\Hh=\begin{pmatrix}q e^{\int_0^t vd\tau} & 0\\
0 &  q^{-1} e^{-\int_0^t vd\tau}\end{pmatrix}.
$$
This gives representation $(c)$ for $h$. Conversely, if $q$, $v$, $h$ satisfy assertions $(a)$-$(c)$, then the Hamiltonian $\Hh = \left(\begin{smallmatrix}h & 0 \\ 0 & 1/h\end{smallmatrix}\right)$ admits $(\mathfrak{q},\mathfrak{v}_1,\mathfrak{v}_2)$ -- factorization $\Hh=G^*QG$ 
for $G$, $Q$ as above. By Theorem \ref{t5}, we have $\Hh \in \hb$. Then Lemma \ref{l42} implies $h \in  A_2(\R_+, \ell^1)$. \qed

\medskip

We now provide some examples of Hamiltonians in class $\hb$. The first two of them show that Theorem~\ref{t1} is essentially sharp.

\medskip

\noindent{\bf Example 1.} Take $\Hh\dd \chi_{[0,L]}\left(\begin{smallmatrix}1&0\\0&0
\end{smallmatrix}\right)+\chi_{(L,\infty)}I_{2\times 2}$, where $L$ is a large integer parameter. Then, $\eta_0=0$, $\eta_j=L+j, j\in \mathbb{N}$ and 
$$
\widetilde \K(\Hh)=\left(\det \left( \int_0^L  \left( \begin{matrix} 1&0\\0&0 \end{matrix}\right)dt+\int_L^{L+2}\left( \begin{matrix} 1&0\\0&1 \end{matrix}\right)dt\right)-4\right)+\sum_{j=1}^\infty (4-4)=2L.
$$
The Titchmarsh-Weyl function can be computed using the formula (2.13) from \cite{BD2017}:
$$
m(z)=\frac{\Phi^+(L,z)+m_L(z)\Phi^-(L,z)}{\Theta^+(L,z)+m_L(z)\Theta^-(L,z)}
,
$$ in which $m_L(z)=i$, because the Titchmarsh-Weyl function of Hamiltonian $I_{2\times 2}$ is equal to constant $i$. Relations \eqref{sd_ff90} yield $m(z)=
\frac{i}{1-izL}.
$
Thus, 
$$
\K_m=-\log(1+L)+\frac{1}{\pi}\int_{\mathbb{R}} \frac{\log(1+x^2L^2)}{1+x^2}dx.
$$
For $L\to\infty$, we can write$$
\frac{1}{\pi}\int_{\mathbb{R}} \frac{\log(1+x^2L^2)}{1+x^2}dx=\frac{1}{\pi}\int_{\mathbb{R}} \frac{\log(x^2L^2)}{1+x^2}dx+O(1)=\frac{2\log L}{\pi}\int_{\mathbb{R}}\frac{dx}{1+x^2}+O(1)=2\log L+O(1).
$$
So, $\K_m= \log L +O(1), \,\, L\to\infty$.
\bigskip

\noindent{\bf Example 2.} Consider Dirac system 
\begin{equation}\label{li170}
JN'(t,z)+V(t)N(t,z)=zN(t,z), \qquad t \in \R_+, \quad z \in \C, \quad N(0,z)=I_{2\times 2},
\end{equation}
with potential
$V=\chi_{[0,T]}\left(\begin{smallmatrix}0 & \eps\\ \eps & 0\end{smallmatrix}\right)$,
where $T$ is a large integer and $\eps$ is a small parameter. They will be chosen such that $L\dd T\eps^2\to\infty$. Define the Hamiltonian 
$\Hh: t \mapsto N^*(t,0)N(t,0)$ on $\R_+$. Then, a straightforward calculation gives
$$
 \Hh(t)=
\begin{pmatrix}
e^{-2\eps t} & 0\\
0 & e^{2\eps t}
\end{pmatrix}, \quad t\in [0,T], \qquad  
\Hh(t)=
\begin{pmatrix}
e^{-2\eps T} & 0\\
0 & e^{2\eps T}
\end{pmatrix}, \quad t\in (T,+\infty).
$$
We have $\eta_n=n$, $n\in \mathbb{Z}_+$ and 
\begin{align*}
\widetilde \K(\Hh)=&\sum_{j=0}^{T-2}\left( \det \left( \int_j^{j+2}
\begin{pmatrix}
e^{-2\eps t} & 0\\
0 & e^{2\eps t}
\end{pmatrix} dt\right)-4\right)+\\
&+\det \left(\int_{T-1}^T 
\begin{pmatrix}
e^{-2\eps t} & 0\\
0 & e^{2\eps t}
\end{pmatrix} dt+\int_{T}^{T+1} 
\begin{pmatrix}
e^{-2\eps T} & 0\\
0 & e^{2\eps T}
\end{pmatrix} dt\right)-4 +\sum_{j=T}^\infty (4-4),\\
=&(T-1)\left(\frac{(1-e^{-4\eps})(e^{4\eps}-1)}{4\eps^2}-4\right) + \left(\left(
\frac{e^{2\eps}-1}{2\eps}+1\right) \left(
\frac{1-e^{-2\eps}}{2\eps}+1\right)   -4\right) \sim 
T\eps^2=L,
\end{align*}
after applying Taylor expansion in small $\eps$. To estimate entropy, we notice that $\Hh$ allows factorization 
$\Hh=G^*QG$ in which $Q=I_{2\times 2}$ and $G=N(t,0)$. Moreover, $V$ is already taken in truncated form similar to \eqref{sd_h77}. The spectral measure $\mu$ of $\Hh$ is absolutely continuous and 
 Lemma \ref{l45} gives
$
\mu'(x)= |\widetilde P^*_{2T}(x)|^{-2},
$
where $(\widetilde P,\widetilde P^*)$ solve Krein system \eqref{li10},\eqref{li11}:
$$
\frac{d}{dr}
\begin{pmatrix}
\widetilde P_{2r}^*\\
\widetilde P_{2r}
\end{pmatrix}
=
\begin{pmatrix} 0&-v\\ -v&2iz \end{pmatrix}
\begin{pmatrix}
\widetilde P_{2r}^*\\
\widetilde P_{2r}
\end{pmatrix}, \quad 
\begin{pmatrix}
\widetilde P_{0}^*\\
\widetilde P_{0}
\end{pmatrix}=
\begin{pmatrix}
 1\\1
\end{pmatrix},
$$
with $v=\eps \chi_{[0,T]}$. We consider $z=x\in [-\eps,\eps]$ and $r\in [0,T]$. Finding eigenvalues $\mu_{\pm}=ix\pm \sqrt{-x^2+\eps^2}$ and eigenvectors $\left(\begin{smallmatrix} \eps\\ -\mu_+
\end{smallmatrix}\right)$,$\left(\begin{smallmatrix} \eps\\ -\mu_-
\end{smallmatrix}\right)$
of matrix 
$\left(\begin{smallmatrix} 0&-\eps\\ -\eps&2ix
\end{smallmatrix}\right)$, we take into account initial data to find
$$
P^*_{2r}(x)=\frac{\eps+\mu_+}{\mu_+-\mu_-}e^{\mu_-r}-\frac{\eps+\mu_-}{\mu_+-\mu_-}e^{\mu_+r},
\quad
P_{2r}(x)=-\frac{\mu_-(\eps+\mu_+)}{\eps(\mu_+-\mu_-)}e^{\mu_-r}+  \frac{\mu_+(\eps+\mu_-)}{\eps(\mu_+-\mu_-)}e^{\mu_+r},
$$
for $r\le T$. We have
$$
\frac{\eps+\mu_+}{\mu_+-\mu_-}e^{\mu_-T}+\frac{\eps+\mu_-}{\mu_--\mu_+}e^{\mu_+T}=\frac{\eps+\mu_-}{\mu_--\mu_+}e^{\mu_+T}\left(1-\frac{\eps+\mu_+}{\eps+\mu_-}e^{(\mu_--\mu_+)T}\right).
$$
Consider $x\in [\frac{\eps}{10},\frac{9\eps}{10}]$. Then, 
$$
\left|\frac{\eps+\mu_-}{\mu_+-\mu_-}\right|\sim 1,\quad
\left|\frac{\eps+\mu_+}{\eps+\mu_-}\right|\sim 1,\quad
-\Re (\mu_-T)\sim \eps T, \quad \Re (\mu_+T)\sim \eps T.
$$
Since $\eps T=(\eps^2T)\eps^{-1}=L\eps^{-1}\to \infty$, this gives us
$$
|P^*_{2T}(x)|^2\sim e^{2\sqrt{\eps^2-x^2}T}
$$
for $x\in [\frac{\eps}{10},\frac{9\eps}{10}]$. Thus, recalling notation \eqref{sd_vi1}, we get  the following estimate
$$
\int_{\frac{\eps}{10}}^{\frac{9\eps}{10}}\frac{\log^- w(x)}{1+x^2}\, dx\sim \eps^2T=L
$$ 
for the spectral measure $\mu = w\,dx + \mus$ of $\Hh$.
We notice that $\lim_{\eps\to 0,T\to\infty}\Hh(t,\eps,T)=\left(\begin{smallmatrix} 1&0\\0&1\end{smallmatrix}\right)$ and this convergence is uniform in $t\in I$ for every fixed segment $I\subseteq \R_+$. Thus,
$\lim_{\eps\to 0,T\to\infty} m(i,\eps,T)= i$. The trivial bound $\log^+a\le a$ yields
$$
\frac{1}{\pi}\int_{\R}\frac{\log^+ w(x)}{1+x^2}\, dx \le 
\frac{1}{\pi}\int_{\R}\frac{w(x)}{1+x^2}\,dx\le
\Im m(i,\eps,T)\lesssim 1,
$$
and thus 
\begin{align*}
\K_m
&=\log \Im m(i,\eps,T)-\frac{1}{\pi}\int_\R \frac{\log w(x)}{1+x^2}\,dx,\\
&\ge-C_1-\frac{1}{\pi}\int_\R \frac{\log^+ \mu' }{1+x^2}dx+\frac{1}{\pi}\int_\R \frac{\log^- w(x) }{1+x^2}\,dx,\\
&\ge -C_2+\frac{1}{\pi}\int_{\frac{\eps}{10}}^{\frac{9\eps}{10}} \frac{\log^- \mu' }{1+x^2}\,dx \sim L,
\end{align*}
for $L\to\infty$. This gives $\K_m\gtrsim L\sim \widetilde \K(\Hh)$. On the other hand, \eqref{ak_1} says that $\K_m\lesssim \widetilde \K(\Hh)$ so  the left-hand side bound in \eqref{ak_1} is sharp up to a constant.\medskip

\noindent {\bf Example 3: Hamiltonians generated by $N_0$.} Consider equation 
\begin{equation}\label{eq92}
J N_0'(t) + V(t)N_0(t) = 0, \quad N_0(0) = I_{2\times 2}, \quad t \in \R_+,
\end{equation} 
from Corollary \ref{t2} in the case when $V=\left(\begin{smallmatrix}v&0\\0&-v\end{smallmatrix}\right)$. Then, we can find $N_0$ and $\Hh=N_0^*N_0$ explicitly. These calculations give
$$
\phi\dd \int_0^t vds,\quad
N_0=\begin{pmatrix}\cosh \phi &\sinh\phi\\  \sinh\phi &\cosh \phi \end{pmatrix}, \quad \Hh=\begin{pmatrix}\cosh (2\phi) &\sinh (2\phi) \\  \sinh (2\phi) &\cosh (2\phi)\end{pmatrix}.
$$
Then, the Theorem \ref{t5} implies  $\Hh\in \hb$ provided that 
$v\in L^2(\R_+)$, because $\Hh$ allows $(\mathfrak{q}, \mathfrak{v}_1, \mathfrak{v}_2)$--factorization in which $G=N_0,Q=I_{2\times 2}$ with $\mathfrak{v}_1=\mathfrak{q}=0$ and $\mathfrak{v}_2\sim \|v\|_2$. 

\medskip

In the case when $V$ takes the form $V=\left(\begin{smallmatrix}0&v\\v&0\end{smallmatrix}\right)$, the equation \eqref{eq92} can also be solved explicitly. That gives yet another class of examples of Hamiltonians in $\hb$. It was discussed in \cite{B2018} in connection with scattering theory for Dirac system.  

\medskip

\noindent {\bf Example 4: Szeg\H{o} condition and indeterminate moment problem.} Consider $\mu$ for which all moments $\{s_k\}$,
\begin{equation}\label{sd_m1}
s_k\dd \int_{\R} x^kd\mu, \quad k\in \mathbb{Z}_+,
\end{equation}
are finite. The sequence $\{s_k\}$ defines the Hamburger moment problem (see \cite{Akh1965} and, e.g., \cite{Berg2014, Berg1995, Berg1995_2},
\cite{Baranov2006}, \cite{Romanov2017} 
for recent developments) and $\mu$ is one of its solutions. We recall that, given a sequence $\{s_k\},k\ge 0$, the moment problem $\{s_k\}\,\rightarrow \,\sigma$ is called indeterminate if, firstly, there is a measure $\sigma$ on the line having $\{s_k\}$ as its moments and, secondly, this measure is not unique. It was noticed by M. Krein, that measures  $\mu$ that satisfy both $\mu\in \sz$ and \eqref{sd_m1} give rise to indeterminate moment problem (see, e.g., \cite{Akh1965}, pp. 87-88). One example of such measures  is $d\mu=w\, dx$, where $w$ is the Freud weight: $w(x)=e^{-|x|^\beta}, \beta>0$, provided that $\beta\in (0,1)$ (see \cite{MK1999} for detailed study of this case). 

\medskip

Every measure that satisfies \eqref{sd_m1} and has  support different from a finite number of points gives rise to a system of polynomials orthogonal on the real line. These polynomials satisfy the three-term recurrence which defines the semi-infinite Jacobi matrix. The inclusion of Jacobi matrices to the more general class of de Branges systems is well-known \cite{Kats1999},\cite{Romanov}. In particular, it shows that measure $\mu\in \sz$ that satisfies \eqref{sd_m1} gives rise to a Hamiltonian $\Hh\in \hb$ for which there is an interval $[0,\ell]$ on which $\rank \Hh=1$. This interval $[0,\ell]$ represents the Jacobi matrix and the elements of $M(t,z)$ in \eqref{eq1} can be expressed in terms of orthogonal polynomials for $t\in [0,\ell)$. However, even for the classical case of Freud weight $w(x)=e^{-|x|^\beta}$ we are not aware of any systematic study of the corresponding Hamiltonian $\Hh(\beta)$ on the interval $[\ell,\infty)$. We notice that our Theorem \ref{t1} yields $\Hh(\beta) \in \hb$ for every $\beta\in (0,1)$.

\medskip

The extensive literature on moment problem contains some cases for which the moments, Jacobi recurrence coefficients, and Nevanlinna matrix of the indeterminate moment problem can be explicitly found. This gives a way of constructing explicit examples of Hamiltonians $\Hh \in \hb$ with known spectral measures in Szeg\H{o} class. For instance, one can consider an example from \cite{Berg_Valent}, Section 2.3, which is related to  birth/death processes. Here, the polynomials $\{F_n\}$ involved satisfy recursion
$$
(\lambda_n+\mu_n-x)F_n=\mu_{n+1}F_{n+1}+\lambda_{n-1}F_{n-1}, \quad F_{-1}=0, \quad F_0=1,
$$
which can be easily symmetrized (see formulas (2.28)--(2.32) in \cite{Berg_Valent}) to produce Jacobi matrix. In the special case when 
$$
\lambda_n=(4n+1)(4n+2)^2(4n+3),\quad \mu_n=(4n-1)(4n)^2(4n+1),
$$
Berg and Valent obtained the asymptotics of $F_n(z)$ for large $n$ and this allowed them to write  (Proposition~3.3.2) the associated Nevanlinna matrix
$$
\left(
\begin{matrix}
A(z) & C(z)\\
B(z) & D(z)
\end{matrix}
\right)
$$
in terms of elementary functions. According to classical theory \cite{Akh1965}, all solutions $\{\mu\}$ to indeterminate moment problem can be parameterized using Nevanlinna matrix in the following way:
$$
\int_{\R} \frac{d\mu_\phi(x)}{z-x}=\frac{A(z)\phi(z)-C(z)}{B(z)\phi(z)-D(z)}, \quad z\in \C_+,
$$
where $\phi$ is arbitrary function from $\N(\C_+)$. In particular, taking $\phi=i$, corresponds to choosing $\Hh=I_{2\times 2}$ on the interval $[\ell,\infty)$, where $\ell$ was mentioned above and can be computed as well. This gives rise to orthogonality measure $\mu_i$ with density determined by the formula (see (2.15) and section 3.5 in \cite{Berg_Valent})
$$
w_i(x)=\frac{1}{\pi(B^2(x)+D^2(x))}.
$$
Since $B$ and $D$ are known, we have (see formula (3.35) in \cite{Berg_Valent})
$$
w_i(x)=
\begin{cases}
C_1(C_2\cos^2 u\cosh^2 u+u^4\sin^2u\sinh^2 u)^{-1}, \quad  &x>0,\\
\widehat C_1 (\widehat C_2(\cos u+\cosh u)^2+u^4(\cos u-\cosh u)^2)^{-1},
&x<0,
\end{cases}
$$
where $u= C_3|x|^{1/4}$ and    $C_1,C_2,\widehat C_1,\widehat C_2,C_3$ are some positive constants known explicitly. Simple analysis shows that $-\log w_i\sim |x|^{1/4}$ which places $\mu_i$ to $\sz$ class.

\section{Reduction to Hamiltonian with unit determinant. Proof of Theorem \ref{l8}}\label{s5}
In this section, we show that the general case in Theorem \ref{t1} can be reduced to the case when Hamiltonian has the unit determinant.
 Our considerations are based on several lemmas that use  additivity  of the entropy function (see assertion $(a)$ in Lemma \ref{l1}) and its upper-semicontinuity. The same ideas were employed in \cite{BD2017}.

\begin{Lem}\label{l6} 
Let $\Hh$, $\Hh_{(k)}$ be singular nontrivial Hamiltonians on $\R_+$ such that $\Hh_{(k)}(t) = \Hh(t)$
for every $k \ge 0$ and all $t \in [0, k]$. Then, we have $\K_{\Hh}(0) \le \limsup_{k\to+\infty}\K_{\Hh_{(k)}}(0)$.
\end{Lem}

\medskip

Lemma \ref{l6} was stated and proved in \cite{BD2017} for diagonal Hamiltonians, see Lemma 4.1 in \cite{BD2017}. Its
proof, however, did not use the fact that the Hamiltonian $\Hh$ is diagonal and hence works in the
general case.

\medskip

\begin{Lem}\label{l10}
The spectral measure of a Hamiltonian $\Hh \in \fc$ lies in $\sz$ class.
\end{Lem}

\beginpf Recall the definition of class $\fc$. By Lemma 2.2, one can assume that $A =\left( \begin{smallmatrix}1&0\\ 0&1\end{smallmatrix} \right)$. Formula (2.14) in \cite{BD2017} then gives
$$
\mu=
\frac{dx}{|F_\ell(x)|^2}, \quad F_{\ell}(z) = \Theta^+ (\ell, z) + i\Theta^-(\ell, z),\quad z \in \mathbb{C},
$$
where $\Theta^{\pm}$ are the entries of the matrix in (2.2). In particular, $F_\ell$ is a function of bounded characteristic
in $\C_+$ and we have $\mu \in \sz$, see Proposition 2.1 in \cite{BD2017}. \qed

\begin{Lem}\label{l7}
Assume that for every Hamiltonian $\Hh\in \fc$ such that $\det \Hh = 1$ a.e.\ on $\R_+$ we have
\begin{equation}
\label{eq15}
c_1\K_\Hh(0)\le \widetilde\K(\Hh)\le c_2\K_\Hh(0)e^{c_2\K_\Hh(0)},
\end{equation}
with an absolute constant c. Then, the same estimates with the same constants $c_1$, $c_2$ hold for every $\Hh \in \fc$.
\end{Lem}
\beginpf Let $\Hh$ be such that $\Hh(t) = A$ for $t \in [\ell, +\infty)$, where $\ell\ge 0$ and $A$ is some  positive 
matrix. For every $\eps > 0$, define $\Hh_{(\eps)}: t \mapsto \Hh(t) + \eps\chi_{[0,\ell]}(t)I_{2\times 2}$ on 
$\R_+$. As before, $I_{2\times 2} = \left(\begin{smallmatrix}  1&0 \\0&1  \end{smallmatrix}\right)$ and $\chi_{[0,\ell]}$ denotes the characteristic function of $[0, \ell].$ For $t > 0$, set
\begin{equation}\label{3.2}
\xi_\eps(t)=\xi_{\Hh_{(\eps)}}(t)=\int_0^t \sqrt{\det \Hh_{(\eps)}(s)} ds,
\end{equation}
and let $\eta_\eps = \eta_{\Hh_{(\eps)}}$ be the inverse function to $\xi_\eps$. Since $\xi_\eps$ bijectively maps $\R_+$ onto $\R_+$, the function $\eta_\eps$ is defined correctly on $\R_+$. Moreover, we have $\det \Hh_{(\eps)} > 0$ a.e.  on $\R_+$, hence $\eta_\eps$ is absolutely continuous on $\R_+$. Consider the Hamiltonian $\widetilde\Hh_{(\eps)}: t\mapsto \eta'_\eps(t)\Hh_{(\eps)}(\eta_\eps(t))$. By construction, $\eta'_\eps(t) = 1/ \sqrt{\det \Hh_{(\eps)}(\eta_\eps(t))}$ a.e.  on $\R_+$, so the Hamiltonian $\widetilde \Hh_{(\eps)}$ has unit determinant a.e.  on $\R_+$. By Lemma~\ref{l10}, the spectral measures $\mu, \mu_{(\eps)}, \widetilde{\mu}_{(\eps)}$
of 
$\Hh$, $\Hh_{(\eps)}$, $\widetilde\Hh_{(\eps)}$, respectively, belong to $\sz$. Our assumption implies the estimates
\begin{equation}\label{eq25}
c_1\K_{\widetilde \Hh_{(\eps)}}(0)\le \widetilde\K( \widetilde \Hh_{(\eps)}  )\le c_2\K_{\widetilde \Hh_{(\eps)}}(0)e^{c_2\K_{\widetilde \Hh_{(\eps)}}(0)}.
\end{equation}
For every $t\ge 0$, we have
$$
\int_{t}^{t+2} \widetilde \Hh_{(\eps)}(s)\,ds=\int_{\eta_\varepsilon(t)}^{\eta_\varepsilon(t+2)}\Hh_{(\eps)}(s)\,ds
$$
by a change of variables. This shows $\widetilde \K (\widetilde \Hh_{(\eps)})=\widetilde \K(\Hh_{(\eps)})$.
Since $\Hh_{(\eps)}$ and $\widetilde\Hh_{(\eps)}$ are equivalent, their Titchmarsh-Weyl functions and spectral measures coincide. Thus, from \eqref{eq25} we get
\begin{equation}\label{eq38}
c_1\K_{\Hh_{(\eps)}}(0)\le \widetilde\K(  \Hh_{(\eps)}  )\le c_2\K_{ \Hh_{(\eps)}}(0)e^{c_2\K_{ \Hh_{(\eps)}}(0)}
\end{equation}
for every $\eps>0$. We claim that $\lim_{\eps\to 0} \widetilde \K(\Hh_{(\eps)})=\widetilde \K(\Hh)$. Let $\eta_n$, $n\ge 0$, be the numbers defined in \eqref{eq2} for the Hamiltonian $\Hh$. Then, $\eta_\eps(n)\leq \eta(n)$ and $\lim_{\eps\to 0} \eta_\eps(n)=\eta_n$ for every $n$. Moreover, for sufficiently large $n_0\ge 0$ we have
$$
\det \int_{\eta_\eps(n)}^{\eta_\eps(n+2)}\Hh_{(\eps)}(t)\, dt=4, \quad  \det \int_{\eta_n}^{\eta_{n+2}}\Hh(t)dt=4,\qquad n\ge n_0, \quad \eps\in (0,1],
$$
due to the fact that $\Hh$, $\Hh_{(\eps)}$ are constant on the corresponding intervals. So, our claim follows from the limiting relations
$$
\lim_{\eps\to 0} \det \int_{\eta_\eps(n)}^{\eta_\eps(n+2)} \Hh_{(\eps)}(s)\, ds=\det \int_{\eta_n}^{\eta_{n+2}} \Hh(s)\, ds, \qquad 0\le n\le n_0,
$$
which are immediate  by the  Lebesgue theorem on  dominated convergence. To complete the proof, it remains to show that $\lim_{\eps\to 0} \K_{\Hh_{(\eps)}}(0)=\K_{\Hh}(0)$. To this end, we will use the following well-known formula (see, e.g., Section 2 in \cite{BD2017}) for $\Hh$ and $\Hh_{(\eps)}$:
\begin{equation}\label{eq10}
m(z)=m_0(z)=\frac{\Phi^+(r,z)+m_r(z)\Phi^-(r,z)}{\Theta^+(r,z)+m_r(z)\Theta^-(r,z)}, \quad r\ge 0, \quad z\in \C_+.
\end{equation}
Denote by $m_{r,\eps}$ the Titchmarsh-Weyl function for $\Hh_{(\eps),r}: t\mapsto \Hh_{(\eps)}(r+t)$. Let also $\Theta^{\pm}_\eps, \Phi^{\pm}_\eps$ be the entries of the solution to Cauchy problem \eqref{eq1} for $\Hh_{(\eps)}$. Since $\Hh_{(\eps)}$ tends to $\Hh$ uniformly on $[0,\ell]$ in the matrix norm and $\Hh=\Hh_{(\eps)}$ on $[\ell,\infty)$, we have
\begin{equation}\label{eq12}
\lim_{\eps\to 0} \Theta^\pm _{\eps}(\ell,i)=\Theta^\pm (\ell,i), \quad \lim_{\eps\to 0} \Phi^\pm _{\eps}(\ell,i)=\Phi^\pm(\ell,i), \quad m_\ell=m_{\ell,\eps} \quad {\rm on} \quad\C_+.
\end{equation}
Applying \eqref{eq10} to $r=\ell$, we get 
$\I_{\Hh}(0)=\Im m_0(i)=\lim_{\eps\to 0} \Im m_{0,\eps}(i)=\lim_{\eps\to 0}\I_{\Hh_{(\eps)}}(0)$. Formula (2.15) in \cite{BD2017} can be rewritten (see also (58) in \cite{B2018}) as
\begin{equation}\label{eq14}
\J_\Hh(0)=\J_\Hh(r)+2\xi_\Hh(r)-2\log |F_r(i)|,
\end{equation}
where $F_r:z\mapsto \Theta^+(r,z)+m_r(z)\Theta^-(r,z)$. The last relation in \eqref{eq12} implies $\J_\Hh(\ell)=\J_{\Hh_{(\eps)}}(\ell)$ while the first two relations together with \eqref{eq14} give us $\J_{\Hh}(0)=\lim_{\eps \to 0} \J_{\Hh_{(\eps)}} (0)$. Recall (see \eqref{sd_h17}) that
$$
\K_\Hh(0) = \log\I_\Hh(0) - \J_\Hh(0), \quad
\K_{\Hh_\eps}(0) = \log\I_{\Hh_\eps}(0) - \J_{\Hh_\eps}(0).
$$
Thus, $\K_{\Hh_{(\eps)}}(0)$ tends to $\K_\Hh(0)$ and the lemma follows from \eqref{eq38}. \qed\medskip

\noindent {\bf Proof  of Theorem \ref{l8}}. By Lemma \ref{l7}, we can drop the condition $\det \Hh=1$ from our assumptions.
Let $\Hh$ be a nontrivial singular Hamiltonian on $\mathbb{R}_+$ such that its spectral measure $\mu$ lies in the class $\sz$. Then, $\K_\Hh(0)< +\infty$. Theorem \ref{p1} and Theorem \ref{p2} imply that $\sqrt{\det \Hh}\notin L^1(\R_+)$. In particular, the sequence $\{\eta_n\}$ in \eqref{eq2} is defined correctly. For $r\ge \eta_2$, consider the Hamiltonian $\widehat \Hh_r\in \fc$, introduced in \eqref{cook59}. The first $[\xi_{\Hh}(r)]-1$ terms defining $\widetilde \K(\Hh)$ and $\widetilde \K(\widehat \Hh_r)$ in \eqref{eq2} are identical. Hence, 
$$
\widetilde \K(\Hh)\le \limsup_{r\to\infty} \widetilde \K(\widehat \Hh_r)\le \limsup_{r\to\infty} c_2\K_{\widehat \Hh_r}(0)e^{c_2\K_{\widehat \Hh_r}(0)}\le c_2\K_{\Hh}(0)e^{c_2\K_{\Hh}(0)},
$$
where the first estimate follows from  construction and definition of $\widetilde \K(\Hh)$, the second inequality follows from assumption of the theorem, and the third one follows from assertion (a) of Lemma~\ref{l1}. Thus, $\Hh\in \hb$ and the first estimate in \eqref{eq15_1} holds.

\medskip

\noindent Conversely, suppose that $\Hh \in\hb$. For every integer $k\ge 0$, define
$$
\widetilde \Hh_{(k)}(t)\dd
\begin{cases}
\Hh(t), & {\rm if}\, t\in [0,\eta_{k+2}],\\
\int_{\eta_{k+2}}^{\eta_{k+3}}\Hh(s)\,ds, & {\rm if}\, t\in (\eta_{k+2}, +\infty).
\end{cases}
$$
For each $t\in \mathbb{Z}_+$, set 
\begin{equation}\label{sd_h21}
\widetilde \eta_{t}=\min\{s: \xi_{\widetilde \Hh_{(k)}}(s)=t\}, 
\end{equation}
where $\xi_{\widetilde \Hh_{(k)}}(s)= \int_0^s \sqrt{\det \widetilde \Hh_{(k)}(\tau)}d\tau$. Then, we have $\widetilde \eta_t=\eta_t$ for every $t\in \{0,\ldots,k+2\}$. By construction, 
\begin{equation}\label{eq51}
\widetilde \K(\widetilde \Hh_{(k)})=\sum_{n=0}^k \left(  
\det \int_{\eta_n}^{\eta_{n+2}}\Hh(s)ds-4 \right)+\det\int_{\widetilde \eta_{k+1}}^{\widetilde\eta_{k+3}}\widetilde \Hh_{(k)}(s)ds-4.
\end{equation}
Indeed, $\widetilde \Hh_{(k)}$ is constant on $[\eta_{k+2},\infty)=[\widetilde\eta_{k+2},\infty)$ and $\Hh=\widetilde\Hh_{(k)}$ on $[0,\eta_{k+2}]$, hence the terms with indices $n\ge k+2$ in formula \eqref{eq2} for $\widetilde\Hh_{(k)}$ vanish, while the terms with indices $n\le k$ coincide with the corresponding terms in \eqref{eq2} for the Hamiltonian $\Hh$.  Since 
$\widetilde\Hh_{(k)}=\int_{\eta_{k+2}}^{\eta_{k+3}}\Hh(s)\,ds$ on $[\eta_{k+2},\infty)$, we have
$$
\int_{\widetilde \eta_{k+1}}^{\widetilde\eta_{k+3}}\widetilde \Hh_{(k)}(s)ds=\left( \int_{\eta_{k+1}}^{\eta_{k+2}}  \Hh(s)ds+(\widetilde\eta_{k+3}-\widetilde\eta_{k+2})\cdot
\int_{\eta_{k+2}}^{\eta_{k+3}}\Hh(s)ds\right)\le \int_{\eta_{k+1}}^{\eta_{k+3}}\Hh(s)\,ds,
 $$
where we used $\widetilde\eta_{k+3}-\widetilde\eta_{k+2}\le 1$.  To obtain the last bound, we recalled \eqref{sd_h21} which gives
$$
(\widetilde\eta_{k+3}-\widetilde\eta_{k+2})
\left(\det \left(\int_{\eta_{k+2}}^{\eta_{k+3}}\Hh(s)ds \right)\right)^{1/2} {=}1,\quad
\widetilde\eta_{k+3}-\widetilde\eta_{k+2}=
\left(\det \left(\int_{\eta_{k+2}}^{\eta_{k+3}}\Hh(s)ds\right) \right)^{-1/2},
$$
and
$$
\det \left(\int_{\eta_{k+2}}^{\eta_{k+3}}\Hh(s)ds\right)\stackrel{\eqref{sd_h10}}{\ge} \left(\int_{\eta_{k+2}}^{\eta_{k+3}}\sqrt{\det \Hh(s)}ds \right)^2\stackrel{\eqref{eq2}}{=}1.
$$
From \eqref{eq51} and  Lemma \ref{sd_h51}, we get $\widetilde \K(\widetilde \Hh_{(k)})\le \widetilde \K(\Hh)$ for every $k\ge 0$. Moreover, 
\begin{equation}\label{sq53}
\lim_{k\to\infty}\widetilde \K(\widetilde \Hh_{(k)})=\widetilde \K(\Hh).
\end{equation}
By Lemma \ref{l10}, the spectral measure of the Hamiltonian $\widetilde \Hh_{(k)}$ belongs to $\sz$ for every $k\ge 0$. Hence, $\mu\in \sz$ and 
$$
\K_\Hh(0)\le\limsup_{k\to\infty}
\K_{\widetilde \Hh_{(k)}} (0)\lesssim
\limsup_{k\to\infty} \widetilde\K( \widetilde \Hh_{(k)})\stackrel{\eqref{sq53}}{\lesssim}
\widetilde \K(\Hh),
$$
where the first inequality follows from Lemma~\ref{l6}. The theorem is proved. \qed

\medskip

\section{Szeg\H{o} condition implies factorization. Proof of Theorem \ref{sd_h32}}\label{s6} 
\noindent In this section we prove Theorem \ref{sd_h32}.
\begin{Lem}\label{sd_h33}
The following estimates are true
\begin{align}
\frac 13 \left| \frac 1x -x   \right|&\le \frac 1x+x-2, \quad x\in (0,1/2)\cup  (2,\infty),\label{li20}\\
\frac x4 &\le \frac 1x+x-2, \quad x\in (0,1/2)\cup  (2,\infty),\label{li21}\\
\frac 29 \left| \frac 1x -x   \right|^2&\le \frac 1x+x-2, \quad x\in [1/2,2] \label{li22}.
\end{align}
\end{Lem}
\beginpf
We will prove the third one, the other bounds can be obtained similarly. Notice that \eqref{li22} is equivalent to showing that
$$
p(x)\dd 2x^4-9x^3+14x^2-9x+2\le 0
$$
for $x\in [1,2]$. We check that $p'(1)=p(1)=p(2)=0$ so factoring gives $p=(2x-1)(x-1)^2(x-2)$ and we get the needed estimate. \qed

\medskip

\noindent {\bf Proof of Theorem \ref{sd_h32}.} 
Recall $\I_{\Hh},\Rr_{\Hh}, \K_{\Hh}$, the functions in $r$,  which were introduced in \eqref{sd_h66} and \eqref{sd_h17}. Let $\Hh = \sth$ and $\det \Hh=1$. Consider
\begin{align*}
G&\dd \left(\begin{matrix}
1/\sqrt{\I_{\Hh}} & \Rr_{\Hh}/\sqrt{\I_{\Hh}}\\
0 &\sqrt{\I_{\Hh}}
\end{matrix}\right),\\
V&\dd \left(\begin{matrix}
0 & \I'_{\Hh}/(2\I_{\Hh})\\
\I'_{\Hh}/(2\I_{\Hh}) & -\Rr'_{\Hh}/\I_{\Hh}
\end{matrix}\right),\\
Q&\dd \left(\begin{matrix}
\I_{\Hh}h_1 & -\Rr_{\Hh}h_1+h   \\
-\Rr_{\Hh}h_1+h & (\Rr^2_{\Hh}h_1-2\Rr_\Hh h+h_2)/\I_\Hh
\end{matrix}\right).
\end{align*}
Now, we can use calculations done in the proof of Lemma 4.3 in \cite{B2018}, to conclude the following:
\begin{itemize}
\item From the last line in (44), \cite{B2018} and  Lemma \ref{sd_h7}, we get $\Hh=G^*QG$. 
\item From the third line  in  (44), \cite{B2018}, we obtain $G'=JVG$.  
\item The fourth  line from the  bottom on the same page gives $\trace Q=2-\K'_\Hh$.
\end{itemize}
Observe that $G\in {\rm SL}(2,\R)$ and $\K_\Hh$ in non-increasing by assertion $(a)$ of Lemma \ref{l1}. Hence, $Q$ is a symmetric matrix with real entries such that 
$$
\trace Q\ge 2, \quad \det Q=\det((G^*)^{-1}\Hh G^{-1})=1,
$$
almost everywhere on $\R_+$. It follows that $Q\ge 0$ almost everywhere on $\R_+$. Moreover, we have
$$
\int_{\R_+} (\trace Q(t)-2) \, dt=-\int_{\R_+} \K'_{\Hh}(t)\,dt=\K_\Hh(0),
$$
where the last equality follows from Lemma \ref{l1}. It remains to estimate the norm of $V$ in $L^1+L^2$. Let $S_1=\{t\in \mathbb{R}_+: 1/2\le \I_\Hh(t)h_1(t)\le 2\}$, $S_2=\R_+\setminus S_1$. Then, we see from the assertion $(d)$ of Lemma \ref{l1} that 
$$
\I'_\Hh/\I_\Hh=\left(\I_\Hh h_1-\frac{1}{\I_\Hh h_1}\right)-\frac{(\Rr'_\Hh/\I_\Hh)^2}{4\I_\Hh h_1}=2g_1+2g_2,
$$
where $2g_1 = (\I_\Hh h_1-\frac{1}{\I_\Hh h_1})\chi_{S_2}-(\Rr'_\Hh/\I_\Hh)^2/(4\I_\Hh h_1)$ and $2g_2 = (\I_\Hh h_1-\frac{1}{\I_\Hh h_1})\chi_{S_1}$. We also define $\widetilde g_1 = -\chi_{S_2}\Rr'_\Hh/\I_\Hh$ and $\widetilde g_2 = -\chi_{S_1}\Rr'_\Hh/\I_\Hh$. Then, we can write $V=V_1+V_2$ with
$$
V_1\dd \left( \begin{matrix} 
0 & g_1\\
g_1 & \widetilde g_1
 \end{matrix}
 \right), 
\quad 
V_2\dd \left( \begin{matrix} 
0 & g_2\\
g_2 & \widetilde g_2
 \end{matrix}
 \right).
$$
Lemma \ref{sd_h33} and assertion $(c)$ of Lemma \ref{l1} imply that
\begin{align*}
2|g_1|&\stackrel{\eqref{li20}}{\le} 3\left(\I_\Hh h_1+\frac{1}{\I_\Hh h_1}-2\right)+\frac{(\Rr'_\Hh/\I_\Hh)^2}{4\I_\Hh h_1} \le -3\K'_\Hh,\\
4|g_2|^2&\stackrel{\eqref{li22}}{\le} \frac 92\left(\I_\Hh h_1+\frac{1}{\I_{\Hh} h_1}-2\right) \le  -9\K'_\Hh/2,\\
|\widetilde g_2|^2&=\chi_{S_1} (\Rr'_\Hh/\I_\Hh)^2\le 8(\Rr'_\Hh/\I_\Hh)^2/(4\I_\Hh h_1) \le -8\K'_\Hh.
\end{align*} 
So, we have three bounds:
$\|g_1\|_{L^1(\R_+)}\le 2\K_\Hh(0), \|g_2\|_{L^2(\R_+)}\le 2\sqrt{\K_\Hh(0)}$, and $\|\widetilde g_2\|_{L^2(\R_+)}\le 3\sqrt{\K_\Hh(0)}$.  Cauchy-Schwarz inequality yields 
$$
\left(\int_{S_2}\left|   \frac{\Rr'_\Hh}{\I_\Hh}\right|\,dt   \right)^2\le 4 \int_{S_2}\left(  \frac{\Rr'_\Hh}{\I_\Hh}\right)^2\frac{dt}{\I_\Hh h_1}\cdot\int_{S_2}\frac{\I_\Hh h_1}{4}\,dt.
$$
By assertion $(c)$ of Lemma \ref{l1} and \eqref{li21}, the right hand side of the above inequality does not exceed $16 \K_\Hh^2(0)$. 
This gives $\|\widetilde g_1\|_{L^1(\R_+)}\le 4 \K_\Hh(0)$. Hence, $V_1\in L^1, V_2\in L^2$ and, moreover,
$$
\|V_1\|_{L^1}\le \|g_1\|_1
\left\|\left(
\begin{matrix}0&1\\1&0\end{matrix}
\right)\right\| +
\|\widetilde g_1\|_1
\left\|\left(
\begin{matrix}0&0\\0&1\end{matrix}
\right)\right\|\le 
6\K_\Hh(0).
$$ Similarly, $\|V_2\|_{L^2}\le  5\sqrt{\K_\Hh(0)}$, as required.
\qed\bigskip

\section{Factorization implies Szeg\H{o} condition. Proof of Theorem \ref{t42}}\label{s7} 
\noindent The  key idea of the proof is to find and estimate an outer function $\widetilde P^*$ defined in $\C_+$, which satisfies
\begin{equation}\label{sd_fa1}
w(x)=|\widetilde P^*(x)|^{-2}
\end{equation} for almost every $x\in \R$. This will provide required bound on the entropy after the multiplicative representation for $\widetilde P^*$ is written at point $i$. We start with some auxiliary statements.
\begin{Lem}\label{l18}
Let $\Hh$ be a singular nontrivial Hamiltonian on $\R_+$ that admits   $(\mathfrak{q}, \mathfrak{v}_1, \mathfrak{v}_2)$ -- factorization.  Then, there exists  $A\in {\rm SL}(2,\R)$ such that
\begin{itemize}
\item[$(a)$]
 $\Hh_A=A^*\Hh A$ admits $(\mathfrak{q}, \mathfrak{v}_1, \mathfrak{v}_2)$ -- factorization as
 $\Hh_A=\widehat G^*\widehat Q\widehat G$, where 
 $\widehat G(0)=
 \left( 
 \begin{smallmatrix} a & 0\\ 0 & a^{-1}
\end{smallmatrix}
\right)$ for some $a\in (0,1]$.
\item[$(b)$]  $m_{A}(i) = i$ for the Titchmarsh-Weyl function $m_A$ of $\Hh_A$.
\end{itemize}
\end{Lem}
\beginpf
Consider $A = G^{-1}(0)BC_\phi$, where $C_\phi=\left( 
 \begin{smallmatrix} \cos \phi & \sin\phi\\ -\sin\phi &\cos\phi
\end{smallmatrix}
\right)$ for some parameter $\phi\in [0,2\pi)$ and $B\in {\rm SL}(2,\R)$  to be chosen later. Define 
$$
\widehat G = C^*_\phi G A, \quad \widehat Q = C^*_\phi Q C_\phi, \quad \widehat V = C^*_\phi V C_\phi.
$$
Using identity $C_\phi^*JC_\phi = J$ (see Lemma \ref{sd_h7}), one can check that
$$
\Hh_A=A^*\Hh A=A^* G^*QGA=\widehat G^*\widehat Q\widehat G, \quad \widehat G'=C_\phi^*JVGA=J\widehat V\widehat G.
$$
Moreover,  $\trace \widehat Q=\trace (QC_\phi C_\phi^*)=\trace Q$ and $\|\widehat V\|_{L^1+L^2}=\|V\|_{L^1+L^2}$. Thus, $\Hh_A$ admits $(\mathfrak{q}, \mathfrak{v}_1, \mathfrak{v}_2)$ -- factorization for any choice of $\phi$ and $B$. Next, choose  symmetric matrix $B\in {\rm SL}(2,\R)$ as 
$$
B =\left( \begin{matrix}
a & b\\
b & d
\end{matrix}\right), \quad d=\frac{I+1}{\sqrt{I(R^2+(I+1)^2)}},\quad a=d\left(I+\frac{R^2}{1+I}\right), \quad b=-\frac{dR}{1+I},
$$
where 
$$
R\dd \Re m_{G^{-1}(0)}(i),  \quad I\dd \Im m_{G^{-1}(0)}(i),
$$
and recall that $m_{G^{-1}(0)}$ denotes Titchmarsh-Weyl function for $\Hh_{G^{-1}(0)}\dd (G^{-1}(0))^*\Hh G^{-1}(0)$. One can verify directly that $\det B=1$.  Then, we  apply \eqref{li16} to check that
$$
m_{\widetilde A}(i)=\frac{dm_{G^{-1}(0)}(i)+b}{bm_{G^{-1}(0)}(i)+a}=i,
$$
where $\widetilde A=G^{-1}(0)B$. Next, we notice that \eqref{li16} implies 
$$
m_{\check\Hh_{C_\phi}}(i)=m_{\check\Hh}(i)=i
$$
for every Hamiltonian $\check \Hh$ and every $\phi$, provided that $m_{\check\Hh}(i)=i$.
Therefore, $m_A(i)=m_{\widetilde A}(i)=i$ for any choice of $\phi\in [0,2\pi)$. From $B=B^*$, $B\in {\rm SL}(2,\R)$, and trace $B\ge 0$, we conclude that $B> 0$. Since $\widehat G(0)=C^*_\phi BC_\phi$, we can take $\phi$ to make sure that $C_\phi$ diagonalizes $B$ and $\widehat G(0)=\left( \begin{smallmatrix} a & 0\\ 0 & a^{-1}\end{smallmatrix}\right)$ for some $a\in (0,1]$. That proves $(a)$ and $(b)$. \qed

\begin{Lem}\label{l12}
Assume that matrix-functions $Q=\left(\begin{smallmatrix}q_1 & q\\ q & q_2   \end{smallmatrix}\right)$, $V=\left(\begin{smallmatrix}v_1 & v\\ v & v_2   \end{smallmatrix}\right)$ satisfy $(a)$--$(d)$ in \eqref{sd_list}. Then, we have $\|v\|_{1,2}+\|v_1\|_{1,2}+\|v_2\|_{1,2}\lesssim  \mathfrak{v}_1+\mathfrak{v}_2$ and $\|q_1-q_2+2iq\|_{1,2}\lesssim \mathfrak{q}+\sqrt{\mathfrak{q}}$.
\end{Lem}
\beginpf We have $\|v\|_{1,2}+\|v_1\|_{1,2}+\|v_2\|_{1,2}\lesssim \|V_1\|_{L^1}+\|V_2\|_{L^2}\le \mathfrak{v}_1+\mathfrak{v}_2$ by definition. Since $\det Q=1$, we also have
\begin{equation}\label{sd_h55}
|q_1-q_2+2iq|^2=(q_1+q_2)^2-4(q_1q_2-q^2)=(q_1+q_2)^2-4.
\end{equation}
Let $S=\{t\in \R_+: q_1(t)\le 3, \, q_2(t)\le 3\}$. At each point of $S$, we have 
$$
|q_1-q_2+2iq|^2=(q_1+q_2-2)(q_1+q_2+2)\le 8(q_1+q_2-2)=8(\trace Q-2).
$$
Therefore, $\|\chi_S(q_1-q_2+2iq)\|_{L^2(\R)}\le 3\sqrt\mathfrak{q}$. We also have $\R_+\setminus S\subseteq \{t:q_1+q_2-2\ge 1\}$ so Chebyshev inequality gives
$$
|\R_+\setminus S|\le \int_{\R_+}(q_1+q_2-2)\,dt=\int_{\R_+}(\trace  Q-2)\,dt\le \mathfrak q.
$$ Using the estimate $|q_1-q_2+2iq|\stackrel{\eqref{sd_h55}}{\le} q_1+q_2=(\trace  Q-2)+2$, we obtain 
$$
\|\chi_{\R_+\setminus S}(q_1-q_2+2iq)\|_{L^1(\R)}\le
\int_{\R_+\setminus S}((\trace  Q-2)+2)\,dt\le \mathfrak q+2|\R_+\setminus S|\le
 3\mathfrak q.
$$ 
The lemma is proved. \qed

\medskip

Next, we will reduce the canonical system with $\Hh$, which admits factorization, to a system of  Dirac type. Then, the system of Dirac type will be further reduced to generalized Krein system. The generalized Krein system turns out to be more convenient for finding representation \eqref{sd_fa1}.

\medskip

Assume that $\Hh$ admits $(\mathfrak{q}, \mathfrak{v}_1, \mathfrak{v}_2)$ -- factorization and $\Hh=G^*QG$. Define $\widetilde \Theta\dd \left(\begin{smallmatrix}\widetilde \Theta^+\\\widetilde \Theta^-\end{smallmatrix}\right)\dd G\Theta$, where $\Theta$ is the first column of the solution to Cauchy problem \eqref{eq1}. Since $JM'=z\Hh M$, we have $(G^*)^{-1}JG^{-1}(GM')=zQ(GM)$. By Lemma \ref{sd_h7}, this yields
$J(GM')=zQ(GM)$, which could be rewritten in the form $J((GM)'-(G'M))=zQ(GM)$. It follows that $J(GM)'+V(GM)=zQ(GM)$, hence
\begin{equation}\label{sd_h45}
J\widetilde \Theta'(t,z)+V(t)\widetilde \Theta(t,z)=zQ(t)\widetilde \Theta(t,z), \quad t \in \R_+, \quad z\in \C,
\end{equation}
for almost every $t\ge 0$. In the case when $Q=I_{2\times 2}$, this equation reduces to Dirac system \eqref{li170}. 

\medskip

Fix absolutely continuous function $u: t\mapsto -\frac 12\int_0^t \trace V(s)\, ds$ on $\R_+$ and consider the following functions for each $r\ge 0$:
\begin{eqnarray}
\widetilde E_r:z\mapsto \widetilde \Theta^+(r,z)+i  \widetilde \Theta^-(r,z),\quad \widetilde P_{2r}:z\mapsto e^{irz-iu(r)}\widetilde E^\sharp_r(z),\label{sd_eqa1} \\
\widetilde E^\sharp_r(z):z\mapsto \widetilde \Theta^+(r,z)-i  \widetilde \Theta^-(r,z), \quad   \widetilde P_{2r}^*:z\mapsto e^{irz+iu(r)}\widetilde E_r(z). \label{eq60}
\end{eqnarray}
\begin{Lem}\label{l16} 
For every $z\in \C$, the function $r\mapsto \widetilde P_{r}^*(z)$ is absolutely continuous in $r$. There are  functions $f(r,z)$ and $g(r)$ that satisfy $$f(\cdot,z)\in L^1(\R_+)+L^2(\R_+), \quad g\ge 0, \quad g\in L^1(\R_+),$$ such that 
\begin{align}
\frac{d}{dr}\widetilde P_{2r}^*(z)&=f(r,z)\widetilde P_{2r}(z)-izg(r)\widetilde P_{2r}^*(z),\label{li10}\\
\frac{d}{dr}\widetilde P_{2r}(z)&=iz(2+g(r)) \widetilde P_{2r}(z)+\overline{f(r,\bar z)}\widetilde P_{2r}^*(z),\label{li11}
\end{align}
for almost every  $r\ge 0$ and all $z\in \C$. Moreover, 
\begin{align}
\|f(\cdot, \pm i)\|_{1,2}&\lesssim \mathfrak{v}_1+\mathfrak{v}_2+\mathfrak{q}+\sqrt{\mathfrak{q}},\\
\|g\|_1&\lesssim  \mathfrak q.\label{sd_f8}
\end{align}
\end{Lem}
\beginpf Define the mapping 
$$
{\bf P}_r(z):r\mapsto \left(\begin{matrix}  \widetilde P_r(z)\\ \widetilde P_r^*(z)\end{matrix}\right), \quad r\ge 0, \quad z\in \C.
$$
We can rewrite \eqref{sd_eqa1} and \eqref{eq60} as
$$
{\bf P}_{2r}(z)=e^{irz}A_1(r)A_2\widetilde \Theta(r,z), \qquad A_1(r)\dd\left( \begin{matrix}
e^{-iu(r)} &0\\
0 & e^{iu(r)}
\end{matrix}\right), \qquad A_2\dd \left(   \begin{matrix} 1 & -i\\
1 & i\end{matrix}\right).
$$
Differentiating with respect to $r$, we get
\begin{align}
{\bf P}_{2r}'
&=iz{\bf P}_{2r}+iu'(r)A_3{\bf P}_{2r}+e^{irz}A_1(r)A_2J^*(zQ(r)-V(r))\widetilde \Theta(r,z)\notag\\
&=(izI_{2\times 2}+iu'(r)A_3+A_1(r)A_2J^*(zQ(r)-V(r))(A_1(r)A_2)^{-1}){\bf P}_{2r},\label{sd_hh1}
\end{align}
where 
$
A_3 \dd \left(\begin{matrix} -1  &0\\0 & 1\end{matrix}\right).
$
Straightforward calculation shows that 
$$
A_1(r)A_2=\left(   \begin{matrix}    e^{-iu(r)} & -ie^{-iu(r)}\\  e^{iu(r)} & ie^{iu(r)}\end{matrix}\right), \quad (A_1(r)A_2)^{-1}=\frac 12\left(\begin{matrix}e^{iu(r)} & e^{-iu(r)}\\  ie^{iu(r)} & -ie^{-iu(r)}\end{matrix}\right),
$$
and
$$
A_1(r)A_2J^*\left(  \begin{matrix} a&c\\c &b\end{matrix}\right)(A_1(r)A_2)^{-1}=\frac 12
\left(
\begin{matrix}
i(a+b)& (2c+ia-ib)e^{-2iu}\\
(2c-ia+ib)e^{2iu} & -i(a+b)
\end{matrix}
\right)
$$
for any $a,b,c\in \C$. Put $a(r)=zq_1(r)-v_1(r), b(r)=zq_2(r)-v_2(r), c(r)=zq(r)-v(r)$, where
$$
Q=\left( \begin{matrix}    q_1 & q\\
q &q_2\end{matrix}\right), \quad V=\left( \begin{matrix}    v_1 & v\\
v &v_2\end{matrix}\right).
$$
Then, $\left(\begin{smallmatrix}a&c\\c&b\end{smallmatrix}\right)=zQ-V$ and \eqref{sd_hh1} shows
\begin{align*}
2\frac{d}{dr}{{ \widetilde P}^*}_{2r}
&=(2c(r)-ia(r)+ib(r))e^{2iu}\widetilde P_{2r}+i(2z+2u'(r)-a(r)-b(r))\widetilde P^*_{2r}\\
&=(2zq-2v+i(v_1-v_2)-iz(q_1-q_2))e^{2iu}\widetilde P_{2r}+i(z(2-q_1-q_2)+2u'+v_1+v_2)\widetilde P_{2r}^*\\
&=(2zq-2v+i(v_1-v_2)-iz(q_1-q_2))e^{2iu}\widetilde P_{2r}+iz(2-\trace  Q(r))\widetilde P_{2r}^*,
\end{align*}
where we used identity $2u'+v_1+v_2=0$. Now, we only need to take
\begin{equation}\label{sd_h63}
f(r,z)\dd zq-v+i(v_1-v_2)/2-iz(q_1-q_2)/2, \quad g(r)\dd \trace Q(r)/2-1,
\end{equation}
to get \eqref{li10}. Formula \eqref{li11} then follows from the relation $$\widetilde P_{2r}(z)=e^{2izr}\overline{\widetilde P^*_{2r}(\bar z)},$$
which can be proved directly by noticing that $\Theta$ and $\widetilde\Theta$ are real for $z\in \R$. Lemma \ref{l12} gives 
$$
\|f(r,z)\|_{1,2}\lesssim  \mathfrak v_1+\mathfrak v_2+\mathfrak q+\sqrt\mathfrak q, \quad z=\pm i,
$$
and we have
$\|g\|_1\lesssim \mathfrak q$ by $(b)$ in \eqref{sd_list}. Function $g$ is non-negative since $\trace Q\ge 2$ (use $\det Q=1$ and $Q>0$ to see this). \qed

\medskip

\noindent {\bf Remark.} Equations \eqref{li10} and \eqref{li11} define the generalization of Krein system. The Krein system was introduced in \cite{Krein54} (see also \cite{Den06}). In fact, \eqref{li10} and \eqref{li11} are identical to Krein system if $g=0$ and $f$ does not depend on $z$.

\medskip

\noindent {\bf Remark.} Consider the dual Hamiltonian $\Hh_d = J^* \Hh J$. Note that if $\Hh$ admits $(\mathfrak q, \mathfrak v_1, \mathfrak v_2)$ factorization $\Hh = G^* Q G$, then the same is true for $\Hh_d = G_d^* Q_d G_d$ with $G_d = J^* G J$, $Q_d = J^* Q J$. This allows us to define the functions $\widetilde P^*_{d,r}$ for $\Hh_d$ as we did it for $\Hh$. The functions $f$, $g$, $f_d$, $g_d$ from the proof of Lemma \ref{l16} for $\Hh$, $\Hh_d$, correspondingly, are related by identities $f_d(r) = - f(r)$, $g_d(r) = g(r)$, $r \ge 0$ due to \eqref{sd_h63} and
$$
Q_d=\begin{pmatrix}    q_2 & -q\\ -q &q_1\end{pmatrix}, \quad V_d\dd J^*VJ=\begin{pmatrix}v_2 & -v\\-v &v_1\end{pmatrix}.
$$
\begin{Lem}\label{l44}
Let $\Hh$ be Hamiltonian which allows $(\mathfrak q, \mathfrak v_1, \mathfrak v_2)$ factorization $\Hh = G^* Q G$. If $G(0) = \left(\begin{smallmatrix}a & 0 \\ 0 &  a^{-1}\end{smallmatrix}\right)$ for some $a>0$, then 
\begin{align}
|\widetilde P^*_r(i)|\le \sqrt 2 ae^{c(\mathfrak q+\mathfrak v_1+\mathfrak v_2^2)}, \quad |\widetilde P^*_{r,d}(i)|\le \sqrt 2 a^{-1}e^{c(\mathfrak q+\mathfrak v_1+\mathfrak v_2^2)}, \label{sd_f11}\\
\sup_{r\ge 0}|\widetilde P^*_r(i) \widetilde P^*_{r,d}(i)|\le 1+ c(\mathfrak{v}_1^2+\mathfrak{v}_2^2+\mathfrak{q})e^{c(\mathfrak{q}+\mathfrak{v}_1+\mathfrak{v}_2^2)}. \label{sd_f12}
\end{align}
\end{Lem}
\beginpf In \eqref{sd_f11}, we will estimate $\widetilde P_r^*$ only, the analysis for $\widetilde P_{r,d}^*$ is analogous. In Lemma \ref{sd_ak3},
take $\Omega$ as
$$
\Omega=
\begin{pmatrix}
g &f(r,i)\\
&\\
\overline{f(r,-i)}& -(2+g)
\end{pmatrix}
$$
and write equations for $(\widetilde P^*,\widetilde P)$ at point $z=i$ in the form:
$$
\frac{d}{dr}\begin{pmatrix}\widetilde P^*\\ \widetilde P\end{pmatrix}=\Omega\begin{pmatrix}\widetilde P^*\\ \widetilde P\end{pmatrix}, \qquad  
\begin{pmatrix}\widetilde P^*_0(i)\\ \widetilde P_0(i)\end{pmatrix}=\begin{pmatrix}a\\a\end{pmatrix}.
$$
Since $g\ge 0$,
$$
\frac{1}{2}(\Omega+\Omega^*)=\left(
\begin{matrix}
g & -v+i(v_1-v_2)/2\\
-v-i(v_1-v_2)/2& -(2+g)
\end{matrix}
\right)\le \left(
\begin{matrix}
g & \mathcal{Q}\\
\overline{\mathcal{Q}}& -2
\end{matrix}
\right),
$$
where $\mathcal{Q}\dd -v+i(v_1-v_2)/2$. Notice that Lemma \ref{l12} gives
\begin{equation}\label{sd_f10}
\mathcal{Q}=\mathcal{Q}^{(1)}+\mathcal{Q}^{(2)},  \quad \|\mathcal{Q}^{(1)}\|_1\lesssim \mathfrak{v}_1, \quad \|\mathcal{Q}^{(2)}\|_2\lesssim \mathfrak{v}_2.
\end{equation}
Let us write $\mathcal{Q}^{(2)}=\mathcal{Q}_1^{(2)}+\mathcal{Q}_2^{(2)}$, where
$$
\mathcal{Q}_1^{(2)}\dd \mathcal{Q}^{(2)}\cdot \chi_{|\mathcal{Q}^{(2)}|>1/10}, \quad \mathcal{Q}_2^{(2)}\dd \mathcal{Q}^{(2)}\cdot \chi_{|\mathcal{Q}^{(2)}|<1/10},
$$
and notice that 
\begin{equation}\label{sd_f9}
\|\mathcal{Q}_1^{(2)}\|_1\le 10\|\mathcal{Q}^{(2)}\|_2^2\lesssim  \mathfrak{v}_2^2.
\end{equation} 
We can write
$$
\frac{1}{2}(\Omega+\Omega^*)\le \left(
\begin{matrix}
g & \mathcal{Q}\\
\overline{\mathcal{Q}}& -2
\end{matrix}
\right)=\Omega_1+\Omega_2\dd \left(
\begin{matrix}
g & \mathcal{Q}_1^{(2)}+\mathcal{Q}^{(1)}\\
\overline{\mathcal{Q}_1^{(2)}}+\overline{\mathcal{Q}^{(1)}}& 0
\end{matrix}
\right)+\left(
\begin{matrix}
0 & \mathcal{Q}_2^{(2)}\\
\overline{ \mathcal{Q}_2^{(2)}  }& -2
\end{matrix}
\right)
$$
and 
$$
\|\Omega_1\|_1\stackrel{\eqref{sd_f8}+\eqref{sd_f9}+\eqref{sd_f10}}\lesssim \mathfrak{q}+\mathfrak{v}_1+\mathfrak{v}_2^2.
$$
The eigenvalues of self-adjoint matrix $\Omega_2$ are $-1\pm\sqrt{1+|\mathcal{Q}_2^{(2)}|^2}$. Since $|\mathcal{Q}_2^{(2)}|<1/10$, we can use Taylor formula to get
$
\Omega_2\lesssim |\mathcal{Q}_2^{(2)}|^2 \cdot I_{2\times 2}.
$
To finish the proof of the first bound in \eqref{sd_f11}, it is left to apply Lemma \ref{sd_ak3}.

Now, consider \eqref{sd_f12}. Denote  $\delta = \mathfrak v_1+\mathfrak v_2+\mathfrak q+\sqrt\mathfrak q$. If $\delta>1$, \eqref{sd_f12} follows from \eqref{sd_f11}. Thus, we can assume that $\delta\le 1$. This implies, in particular, that $\max\{\mathfrak{q},\mathfrak{v}_1,\mathfrak{v}_2\}\le 1$ and we only need to show that
\begin{equation}\label{sd_f15}
\sup_{r\ge 0}|\widetilde P^*_r(i) \widetilde P^*_{r,d}(i)|\le 1+ c\delta^2.
\end{equation}
If $f(r,z)$ is the function from Lemma \ref{l16}, we let $f(r) = f(r,i)$ and $\widetilde f(r) = \overline{f(r,-i)}$ for all $r\ge 0$. Define 
$$
\kappa(r)= \int_0^r g(t)dt,\quad 
p^*(r)= a^{-1}e^{-\kappa(r)}\widetilde P^*_{2r}(i), \quad p(r)= a^{-1}e^{2r+\kappa(r)}\widetilde P_{2r}(i), \quad r\ge 0,
$$
where $g$ was introduced in \eqref{sd_h63}.
Then, we have $p^*(0) = p(0) = 1$, ${p^*}'(r)=e^{-2r-2\kappa}fp$ and $p'=e^{2r+2\kappa}\widetilde fp^*$ a.e.  on $\R_+$. It follows that
\begin{align}
p^*(r)&=1+\int_0^r e^{-2t-2\kappa(t)}f(t)p(t)\,dt,\notag\\
&=1+\int_0^r e^{-2t-2\kappa(t)}f(t)\left(1+\int_0^te^{2s+2\kappa(s)} \widetilde f(s)p^*(s)ds\right)\,dt, \label{li14}  \\
&=1+\int_0^r e^{-2t-2\kappa(t)}f(t)dt+\int_0^r \widetilde f(s)p^*(s)\int_s^r f(t)e^{2(s-t+\kappa(s)-\kappa(t))}\,dt\,ds.\notag
\end{align}
Using $g\ge 0$, we obtain $\kappa(t)\ge 0$ and $\kappa(s)-\kappa(t)\le 0$, so
$$
|p^*(r)|\le 1+\int_0^r e^{-2t}|f(t)|\,dt+\int_0^r |\widetilde f(s)p^*(s)|\int_s^\infty |f(t)|e^{2s-2t}dt\,ds.
$$
Now we can apply Gr\"onwall inequality to get
\begin{align*}
\sup_{r\ge 0}|p^*(r)|
&\le \left(1+\int_0^\infty e^{-2t}|f(t)|\,dt\right)\exp \left(\int_0^\infty \int_s^\infty |\widetilde f(s)f(t)|e^{2s-2t}\,dt\,ds\right),\\
&\le (1+\|f\|_{1,2})\exp\left(\int_{\R}\int_{\R} |\widetilde f(s)|\cdot \chi_{\R_+}(t-s)e^{-2(t-s)}\cdot|f(t)|\,dt\,ds \right),
\end{align*}
where we extended $f,\widetilde f$ to $(-\infty, 0)$ by zero. From Young's inequality for convolution, i.e.,
$$
\int_\R \int_\R h_1(s)h_2(t-s)h_3(t)dsdt\le \|h_1\|_{p_1}\|h_2\|_{p_2}\|h_3\|_{p_3}, \quad 
p_1^{-1}+p_2^{-1}+p_3^{-1}=2, \quad p_j\ge 1,\quad j=1,2,3,
$$
we obtain
$$
 \int_{\R}\int_{\R} |\widetilde f (s)|    \cdot \chi_{\R_+}   (t-s)  e^{-2(t-s)} \cdot |f(t)|dtds\le \|f\|_{1,2}\|\widetilde f\|_{1,2}.
$$
It follows that $\sup_{r\ge 0}|p^*(r)|\le (1+\|f\|_{1,2})e^{\|f\|_{1,2}\|\widetilde f\|_{1,2}}$. By Lemma \ref{l16}, we have 
$\|f\|_{1,2} \lesssim  \delta$, hence $\sup_{r\ge 0}|p^*(r)|\le (1+c\delta)e^{c\delta^2}$ with an absolute constant $c$. 
The same argument applies to the ``dual'' function $p_d^*(r)\dd ae^{-\kappa(r)}\widetilde P^*_{2r,d}(i)$. In particular, we have 
$$
p_d^*(r) = 1-\int_0^r e^{-2t-2\kappa(t)}f(t)dt+\int_0^r \widetilde f(s)p_d^*(s)\int_s^r f(t)e^{2(s-t+\kappa(s)-\kappa(t))}dtds, 
$$
where we used formula \eqref{li14} and the fact that $f_d = - f$ (see Remark before the proof). It follows that $\sup_{r\ge 0}|p_d^*(r)|\le (1+c\delta)e^{c\delta^2}$. Multiplying $p^*$ with $p_d^*$, we see that the linear in $f$ terms cancel out, while the other terms are controlled by $\delta$. This yields the following estimate:    
\begin{equation}\label{li15}
|p^*(r)p_d^*(r)| \le  1+C\delta^2e^{C\delta^2}
\end{equation}
after combining all terms.
Since $\int_{0}^{+\infty}g(t)\,dt \lesssim  \mathfrak q$ and $p^*(r)p_d^*(r) = e^{-2\kappa(r)}\widetilde P^*_{2r}(i)\widetilde P^*_{2r,d}(i)$, we see that \eqref{li15} implies \eqref{sd_f15}, because
$$
e^{C\mathfrak q}\le 1+C\mathfrak qe^{C\mathfrak q}, \quad (1+C\mathfrak qe^{C\mathfrak q})(1+C\delta^2e^{C\delta^2})\le 1+c\delta^2,
$$ 
for $\delta\le 1$. \qed

\medskip

Given $\ell\ge 0$ and a Hamiltonian $\Hh$ which allows $(\mathfrak{q}, \mathfrak{v}_1, \mathfrak{v}_2)$ -- factorization $\Hh=G^*QG$, we can define the following approximation (compare it with \eqref{cook59} which always exists):
\begin{equation}\label{sd_h77}
\Hh_{(\ell)}\dd (G_{(\ell)})^*Q_{(\ell)}G_{(\ell)},
\end{equation}
where $Q_{(\ell)}=Q\chi_{[0,\ell]}+I_{2\times 2}\chi_{[\ell,\infty)}$ and $G_{(\ell)}$ solves Cauchy problem
$$
G_{(\ell)}'=JV_{(\ell)}G_{(\ell)}, \quad G_{(\ell)}(0)=G(0),
$$
for $V_{(\ell)}\dd V\chi_{[0,\ell]}$. From the definition, we get $\Hh(t)=\Hh_{(\ell)}(t)$ for $t\in [0,\ell]$ and 
$\Hh_{(\ell)}(t)=G^*(\ell)G(\ell)$ for $t\in [\ell,\infty)$. Clearly, Hamiltonian $\Hh_{(\ell)}$ admits  $(\mathfrak{q}, \mathfrak{v}_1, \mathfrak{v}_2)$--factorization.
Moreover, \eqref{sd_h63} shows that $f_{(\ell)}(t)=g_{(\ell)}(t)=0$ for $t>\ell$. Therefore, \eqref{li10} yields $\widetilde{P}^*_{2r,(\ell)}=\widetilde P^*_{2\ell}$ for $r>\ell$. In the next lemma, we show that $1/{{\widetilde{P}}^*_{2\ell}}$ is the function we are looking for: an outer function in $\C_+$ which provides a factorization of the spectral measure of $\Hh_{(\ell)}$.

\begin{Lem}\label{l45}
Let $\Hh$ be a Hamiltonian which allows $(\mathfrak{q}, \mathfrak{v}_1, \mathfrak{v}_2)$-factorization. Let $\widetilde P^*_{2\ell}$ be defined by \eqref{eq60} for $r=\ell$. Then, $\widetilde P^*_{2\ell}$ satisfies the following properties:
\begin{itemize}
\item[$(a)$] $|\widetilde P^*_{2\ell}(x)|^{-2}\,dx$ is the spectral measure for $\Hh_{(\ell)}$,
\item[$(b)$] $\widetilde P_{2\ell}^*(z)$ is an outer function in $z\in \C_+$ and, in particular,
\begin{equation}\label{eq20}
\frac 1\pi \int_{\R}\log|\widetilde P^*_{2\ell}(x)|^2\frac{\Im z}{|x-z|^2}dx=\log|\widetilde P_{2\ell}^*(z)|^2,\quad z\in \C_+.
\end{equation}
\end{itemize}
\end{Lem}
\beginpf If
$G_{(\ell)}=\begin{pmatrix} g_{11} & g_{12}\\ g_{21} & g_{22} \end{pmatrix}$,
then \eqref{sd_h62} and \eqref{li16} imply that the Titchmarsh-Weyl function for Hamiltonian $G_{(\ell)}^*G_{(\ell)}$ is given by $\displaystyle\frac{g_{22}i+g_{12}}{g_{21}i+g_{11}}$, since the Titchmarsh-Weyl function of Hamiltonian $I_{2\times 2}$ is equal to $i$. Therefore, the density of the spectral measure of $\Hh_{(\ell)}$ can be written as follows (see $(b)$, Lemma 2.2 in \cite{B2018}):
\begin{align*}
w_{(\ell)}(x)&=\frac{1/|g_{21}i+g_{11}|^2}{\left|\Theta^+(\ell,x)+\Theta^-(\ell,x)\frac{g_{22}i+g_{12}}{g_{21}i+g_{11}}\right|^2}\\
&=\frac{1}{|(\Theta^+(\ell,x)g_{11}+\Theta^-(\ell,x)g_{12})+i(\Theta^+(\ell,x)g_{21}+\Theta^-(\ell,x)g_{22})|^2}\\
&=\frac{1}{|\widetilde\Theta^+(\ell,x)+i\widetilde\Theta^-(\ell,x)|^2}\stackrel{\eqref{eq60}}{=}\frac{1}{|\widetilde P^*_{2\ell}(x)|^2},
\end{align*}
which proves $(a)$. Recall that $\widetilde \Theta=G\Theta=\left(   \begin{smallmatrix}  \widetilde\Theta^+\\
\widetilde\Theta^-\end{smallmatrix}\right)$. By definition of $J$, we have 
$$
2i\Im (\widetilde \Theta^+(r,z)\overline{\widetilde \Theta^-(r,z)})=\langle J\widetilde \Theta (r,z),\widetilde \Theta(r,z)\rangle_{\C^2}.
$$ 
Since $G\in {\rm SL}(2,\R)$, we can apply Lemma \ref{sd_h7} to get 
$$
\langle J\widetilde \Theta (r,z),\widetilde \Theta(r,z)\rangle_{\C^2}=\langle J\Theta (r,z),\Theta(r,z)\rangle_{\C^2}.
$$
It follows that
\begin{equation}\label{sd_h52}
\Im (\widetilde \Theta^+(r,z)\overline{\widetilde \Theta^-(r,z)}) = \Im (\Theta^+(r,z) \overline{\Theta^-(r,z)}), \qquad z \in \C_+.
\end{equation}
Let $E_r(z) = \Theta^+(r,z)+i\Theta^-(r,z)$ and notice that
\begin{align}
|\widetilde E_r(z)|^2
&=|\widetilde \Theta^+(r,z)+i\widetilde \Theta^-(r,z)|^2=\|G(r)\Theta(r,z)\|^2_{\C^2}+2\Im (\widetilde \Theta^+(r,z)\overline{\widetilde\Theta^-(r,z)}) \notag\\
&=\|G(r)\Theta(r,z)\|^2_{\C^2}+2\Im ( \Theta^+(r,z)\overline{\Theta^-(r,z)}),  \label{sd_h53}\\
|E_r(z)|^2
&=|\Theta^+(r,z)+i \Theta^-(r,z)|^2=\|\Theta(r,z)\|^2_{\C^2}+2\Im (\Theta^+(r,z)\overline{\Theta^-(r,z)}).\notag
\end{align} 
Since $G(r)\in {\rm SL}(2,\R)$, it is invertible and we have
$$
0<c_{(G(r))}\|\Theta(r,z)\|^2_{\C^2}\le {\|G(r)\Theta(r,z)\|^2_{\C^2}}\le C_{(G(r))}\|\Theta(r,z)\|^2_{\C^2}
$$
for all $z\in \C$. Formulas \eqref{sd_h53} then yield
$$
0<\widetilde c_{(G(r))}|E_r(z)|\le {|\widetilde E_r(z)|}\le \widetilde C_{(G(r))}|E_r(z)|.
$$
On the other hand, it is known that the entire function $E_r$ is in Hermite-Biehler class. In particular, it has no zeroes in $\C_+$, which implies that $\widetilde E_r$ and $\widetilde P^*_{r}$ have no zeroes in $\C_+$ as well. It is also known (see, e.g., Section 6 in \cite{Romanov}) that $E_r$ has bounded type in $\C_{+}$ and, moreover, 
$$
\limsup_{y\to +\infty}\frac{\log|E_r(iy)|}{y}=\int_0^r \sqrt{\det \Hh(t)}dt=r.
$$
Therefore, the function $\widetilde P_{2r}^*=e^{irz+iu(r)}\widetilde E_r(z)$ has bounded type in $\C_+$ as well (in particular, $\log |\widetilde P_{2r}^*(x)|\cdot (x^2+1)^{-1}\in L^1(\R)$), and
\begin{equation}\label{eq19}
\limsup_{y\to +\infty}\frac{\log|\widetilde P^*_{2r}(iy)|}{y}=\limsup_{y\to +\infty}\frac{\log|e^{-ry}E_r(iy)|}{y}=0.
\end{equation}
Since $\widetilde P^*_{2r}$ is of bounded type, it allows canonical factorization (see Theorem 5.5 in \cite{Garnett}):
$$
\widetilde P^*_{2r}(z)=\frac{B_1(z)I_1(z)}{B_2(z)I_2(z)}O(z),\quad z\in \C_+,
$$
where $B_{1},B_{2}$ are  Blaschke products, $I_{1},I_2$ are inner functions, and $O$ is an outer function. However, since $\widetilde P^*_{2r}$has no zeroes and it is entire, we get $B_1=B_2=1$ and $I_1/I_2=e^{i\beta z},\beta\in \R$. Then, \eqref{eq19} shows that $\beta=0$ so 
$\widetilde P_{2r}^*$ is an outer function in $\C_+$ and \eqref{eq20} holds. \qed

\medskip

\noindent {\bf Proof of Theorem \ref{t42}.} 
By Lemma \ref{l6}, it is sufficient to consider $\Hh_{(\ell)}$ and prove $$ 
\K_{\Hh_{(\ell)} }(0)\lesssim \min\{\mathfrak{v}_1,\mathfrak{v}_1^2\}+\mathfrak{v}_2^2+\mathfrak{q}
$$ 
for all $\ell$. Denote Titchmarsh-Weyl function of $\Hh_{(\ell)}$ by $m_{(\ell)}$. By  Lemma~\ref{l18}, we may additionally assume that $m_{(\ell)}(i)=i$ and that $\Hh_{(\ell)}$ admits $(\mathfrak{q}, \mathfrak{v}_1, \mathfrak{v}_2)$ -- factorization $\Hh_{(\ell)}=G^*QG$ with $G(0)=\left(\begin{smallmatrix} a & 0\\ 0 & a^{-1}\end{smallmatrix}\right)$ for some $a>0$. We notice here that if Hamiltonian is an approximation of the type \eqref{sd_h77}, it will be of the same type after modifying it as in Lemma~\ref{l18}. Using  Lemma \ref{l45}, we obtain 
$$
\K_{\Hh_{(\ell)}}(0)\stackrel{\eqref{sd_h17}}{=}\log \I_{\Hh_{(\ell)}}(0)- \J_{\Hh_{(\ell)}}(0)=
-\frac 1\pi \int _\R \frac{\log |\widetilde P^*_{2\ell}(x)|^{-2}}{x^2+1}dx=\log |\widetilde P^*_{2\ell}(i)|^2. 
$$
By Corollary \ref{c1}, we have $\I_{\Hh_{(\ell,d)}}(0) = \Im (-1/i) = 1$ and $\K_{\Hh_{(\ell,d)}}(0) = \K_{\Hh_{(\ell)}}(0)$ for the dual Hamiltonian ${\Hh_{(\ell,d)}}$. Hence, $\K_{\Hh_{(\ell,d)}}(0) = \log |\widetilde P^*_{2\ell,d}(i)|^2$. Then, Lemma \ref{l44}  gives the estimate
\begin{equation}\nonumber
\K_{\Hh_{(\ell)}}(0) = \frac{\K_{\Hh_{(\ell)}}(0) + \K_{\Hh_{(\ell,d)}}(0)}{2} = \log |\widetilde P^*_{2\ell}(i) \widetilde P^*_{2\ell,d}(i)| \lesssim 
\left\{
\begin{array}{cc}
\mathfrak{v}_1^2+\mathfrak{v}_2^2+\mathfrak{q}, &{\rm if}\,\max\{\mathfrak{v}_1,\mathfrak{v}_2,\mathfrak{q}\}\le 1,\\
\mathfrak{v}_1+\mathfrak{v}_2^2+\mathfrak{q}, &{\rm if}\,\max\{\mathfrak{v}_1,\mathfrak{v}_2,\mathfrak{q}\}>1
\end{array} \right.
\end{equation}
and the theorem follows. 
 \qed\smallskip
 
 \noindent {\bf Remark.} In this paper, we do not develop the full Szeg\H{o} theory for generalized Krein systems. In particular, we do not study convergence of $\widetilde P_r^*(z)$ to Szeg\H{o} function when $z\in \C_+$. In \cite{B2018}, this question was addressed for a special kind of Krein systems. We believe that the same argument  works in full generality.

\section{Factorization controls mean oscillation. Proof of Theorem \ref{t5}}\label{s8}
In this section, we show that a Hamiltonian which admits $(\mathfrak{q}, \mathfrak{v}_1, \mathfrak{v}_2)$--factorization belongs to $ \hb$.

\medskip 

\noindent{\bf Proof of Theorem \ref{t5}.} Suppose that $\Hh = G^* Q G$ is  $(\mathfrak{q}, \mathfrak{v}_1, \mathfrak{v}_2)$--factorization of $\Hh$. Take $n \in \Z_+$ and define $U_n(t)$ as the solution to
$$
U_n'=JVU_n,\quad U_n(n)=I_{2\times 2}.
$$
Then, we have $G(t) = U_n(t)G(n)$ for $t\ge n$.
Defining 
$$
W_n(t) =  \int_{n}^{t}JV(s)\,ds, \qquad \Delta^{(1)}_n(t) = \int_n^t JV(s_1)\int_n^{s_1}JV(s_2)U_n(s_2)\,ds_2\,ds_1,
$$
we iterate integral equation
$$
U_n(t)=I_{2\times 2}+\int_n^t JV(s)U_n(s)ds
$$
once to write $G$ in the form
$$
G(t) = \Bigl(I_{2\times 2} + W_n(t) + \Delta_{n}^{(1)}(t)\Bigr)G(n), \qquad t\ge n.
$$
Since $G(n)\in {\rm SL}(2,\R)$, we get
\begin{align*}
\det \int_{n}^{n+2}\Hh(t)\,dt 
&= \det \left(\int_n^{n+2}G^*(n)\Bigl(I_{2\times 2}+W_n^*(t)+{\Delta^{(1)}_n(t)}^*\Bigr)Q(t)\Bigl(I_{2\times 2}+W_n(t)+\Delta_n^{(1)}(t)\Bigr)G(n)\,dt\right)\\
&=\det \left(\int_n^{n+2}\Bigl(I_{2\times 2}+W_n^*(t)+{\Delta^{(1)}_n(t)}^*\Bigr)Q(t)\Bigl(I_{2\times 2}+W_n(t)+\Delta_n^{(1)}(t)\Bigr)\,dt\right).
\end{align*}
Denote $\mathfrak v = \mathfrak v_1 + \mathfrak v_2$. Since $\|V\| \le \|V_1\|+\|V_2\|$, we have
$$
\int_n^{n+2}\|V\|dt\le \int_n^{n+2}(\|V_1\|+\|V_2\|)dt \le \|V_1\|_1+\sqrt 2\|V_2\|_2\le 2\mathfrak{v},
$$
by Cauchy-Schwarz inequality. It follows that 
\begin{align}
\sup_{t\in [n,n+2]}\|U_n(t)\| &\le \exp\left(\int_n^{n+2}\|V(t)\|dt\right)\lesssim e^{2\mathfrak{v}}, \quad
\sum_{n \ge 0}\sup_{t\in [n,n+2]}\|\Delta_n^{(1)}(t)\| \lesssim \mathfrak{v}^2e^{2\mathfrak{v}}, \label{li18}\\
\int_{n}^{n+2}\|Q(t)\|\,dt &\le \int_{n}^{n+2}\trace Q(t)\,dt \lesssim 1+\mathfrak{q}, \quad 
\sup_{t \in [n,n+2]}\|W_n(t)\| \le \int_{n}^{n+2}\|V(t)\|\,dt \lesssim \mathfrak{v}. \label{li19}
\end{align}
For $2 \times 2$ matrices $A$ and $B$, we have 
\begin{equation}\label{sd_h43}
\det (A+B)= \det A+O(\|B\|^2+\|A\|\cdot \|B\|),\end{equation}
as can be easily seen by doing calculation on the  determinant. So, taking
$$
A_n = \int_n^{n+2}\bigl(I_{2\times 2}+W_n^*\bigr)Q\bigl(I_{2\times 2}+W_n\bigr)\,dt, \qquad  B_n = \int_n^{n+2}\bigl(U_n^*Q{\Delta_n^{(1)}}+{\Delta^{(1)}_n}^*QU_n\bigr)\,dt,
$$
we get $$\det \int_{n}^{n+2}\Hh(t)\,dt = \det(A_n+B_n) \dd  \det A_n + \delta_n.$$ Notice, that the sum of smaller order terms $\delta_n$ allows the bound
\begin{equation}\label{sd_h57}
\sum_{n\ge 0}|\delta_n| \lesssim \mathfrak{v}^2(1+\mathfrak{q})^2e^{9\mathfrak{v}},
\end{equation}
as follows from \eqref{li18}, \eqref{li19}, and \eqref{sd_h43}. Since
$$
\sum_{n \ge 0}\left(\det \int_{n}^{n+2} \Hh(t)\,dt - 4 \right) \le \sum_{n \ge 0}\left(\det A_n - 4 \right) + \sum_{n\ge 0}|\delta_n|,
$$
it remains to estimate $\det A_n$. We have
\begin{align*}
\frac 12A_n&=I_{2\times 2}+\frac 12\left(   \int_n^{n+2}\Bigl(W_n^*(t) + W_n(t)\Bigr)\,dt  +\int_{n}^{n+2}(Q(t)-I_{2\times 2})\, dt \right)+\Delta_n^{(2)},\\
\Delta_n^{(2)}&\dd \frac 12\int_n^{n+2}\Bigl(W_n^*(t)Q(t)W_n(t)+(Q(t)-I_{2\times 2})W_n(t)+W_n^*(t)(Q(t) - I_{2\times 2})\Bigr)\,dt.
\end{align*}
Let $\lambda(t)$, $\lambda^{-1}(t)$ denote the eigenvalues of the matrix $Q(t)$ for $t \ge 0$ and we can assume that $\lambda(t) \ge 1$, because $\det Q=1$. Then $\lambda(t) + \lambda^{-1}(t) - 2=\trace  Q(t)-2$ is a non-negative function whose integral over $\R_+$ does not exceed $\mathfrak q$. Define the function $k: t \mapsto |\lambda(t) - 1| + |\lambda^{-1}(t) - 1| = \lambda(t) - \lambda^{-1}(t)$ on $\R_+$ and observe that $\|Q(t) - I_{2\times 2}\|=\max\{|\lambda(t)-1|,|\lambda^{-1}(t)-1|\} \le k(t)$. Consider two sets $S = \{t\in \R_+: \lambda(t) \ge 2\}$ and $S^c = \R_+ \setminus S$. We define $k_1 = \chi_S k$, $k_2 = \chi_{S^c} k$ and use \eqref{li20} and \eqref{li22} to get
$$
\|k_1(t)\|_{L^1(\R)} \stackrel{\eqref{li20}}{\le} 3\mathfrak q, \qquad \|k_2(t)\|_{L^2(\R)}^{2} \stackrel{\eqref{li22}}{\le} \frac{9\mathfrak q}{2}, \qquad \|k_2(t)\|_{L^{\infty}(\R)} \le \frac 32.
$$
Recall that $V=V_1+V_2$ and introduce
\begin{align*}
c_{1, n} =  \int_{n}^{n+2}\|V_1(t)\|\,dt&, \qquad c_{2, n} = \sqrt{\int_{n}^{n+2}\|V_2(t)\|^2\,dt},\\
d_{1,n} = \int_{n}^{n+2} k_1(t)\,dt&, \qquad d_{2,n} = \sqrt{\int_{n}^{n+2} k^2_2(t)\,dt}.
\end{align*}
Then, Cauchy-Schwarz inequality gives
\begin{align*}
\sup_{t \in [n, n+2]}\|W_n(t)\| \lesssim c_{1,n} + c_{2,n}&, \qquad \int_{n}^{n+2}\|Q(t) - I_{2\times 2}\|dt \lesssim d_{1,n} + d_{2,n},\\
\sum_{n\ge 0} c_{1, n} \lesssim \mathfrak{v}, \qquad \sum_{n\ge 0} c_{2,n}^2 \lesssim \mathfrak{v}^2 &, \qquad \sum_{n\ge 0} d_{1, n} \lesssim \mathfrak{q}, \qquad \sum_{n\ge 0} d_{2,n}^2 \lesssim \mathfrak{q}.
\end{align*}
Hence,
\begin{align*}
\|\Delta_n^{(2)}\|&\lesssim (1 + \mathfrak{q}) (c_{1,n} + c_{2,n})^2 +  (c_{1,n} + c_{2,n})(d_{1,n} + d_{2,n}),\\ 
\sum_{n\ge 0} \|\Delta_n^{(2)}\| &\lesssim (1 + \mathfrak{q})\sum_{n \ge 0} (c_{1,n}^{2} + c_{2,n}^{2}) 
+\sum_{n\ge 0}(c_{1,n} + c_{2,n})(d_{1,n} + d_{2,n}) 
\lesssim (1 + \mathfrak{q})\mathfrak{v}^2 +(\mathfrak{q}+\mathfrak{q}^{1/2})\mathfrak{v}.
\end{align*}
Notice that $$\det (I_{2\times 2}+B)=1+\trace  B+O(\|B\|^2)$$ for any $2\times 2$ matrix $B$ and
$$\trace(VJ^*+JV)=\trace(V(J+J^*))=0.$$
So, we are left with
\begin{align*}
\det (A_n/2)
&=1+\frac 12 \int_{n}^{n+2}\bigl(\trace Q(t)-2\bigr)\,dt + \trace \Delta_n^{(2)} + O\left(c_{1,n}^2+c_{2,n}^2+d^2_{1,n} + d^2_{2,n}\right)+O(\|\Delta_n^{(2)}\|^2)\,\\
&\dd 1+I_1+I_2+I_3+I_4.
\end{align*}
Since 
$
\trace \Delta_n^{(2)} \le 2\|\Delta_n^{(2)}\|$, we have
$$
\sum_{n\ge 0}(\det A_n - 4) \lesssim \overbrace{\mathfrak{q}}^{\sum_n I_1} + \overbrace{\Bigl((1+\mathfrak{q})\mathfrak{v}^2+(\mathfrak{q}+\mathfrak{q}^{1/2})\mathfrak{v}\Bigr)}^{\sum_n I_2}+\overbrace{\Bigl(  \mathfrak{v}^2+\mathfrak{q}+\mathfrak{q}^2 \Bigr)}^{\sum_n I_3}+\overbrace{\Bigl((1+\mathfrak{q})^2\mathfrak{v}^4+\mathfrak{v}^2(\mathfrak{q}+\mathfrak{q}^2)\Bigr)}^{\sum_n I_4}.
$$
Combining it with \eqref{sd_h57} and using a trivial bound $2\mathfrak{v}\mathfrak{q}^{1/2}\le (\mathfrak{v}^2+\mathfrak{q})$, we get
$\widetilde{\mathcal{K}}_\Hh\le c(\mathfrak{q}+\mathfrak{q}^2+\mathfrak{v}^2)e^{c\mathfrak{v}}$ with an absolute constant $c$. \qed

\bigskip

\section{The condition on mean oscillation implies factorization. Proof of Theorem \ref{t6}}\label{s9}
Now, we turn to proving Theorem \ref{t6}. We need first one auxiliary result on triangular factorization of matrices. Suppose $A= \left(\begin{smallmatrix}  a_1&a\\a&a_2\end{smallmatrix}\right)$ is positive real $2\times 2$ matrix. We denote by $\Lambda_A$ real upper-triangular matrix which satisfies
\begin{equation}\label{sd_f7}
\Lambda_A=\left(\begin{smallmatrix} \lambda_1(a)& \lambda(a)\\0&\lambda_2(a) \end{smallmatrix}\right),\quad \lambda_1(a)>0, \quad \lambda_2(a)>0, \quad \Lambda_A^*\Lambda_A=A.
\end{equation}
One can solve equations for $\lambda,\lambda_1,\lambda_2$ and find $\Lambda_A$ uniquely:
$$
\Lambda_A=\left(
\begin{matrix}
\sqrt{a_1} & a/\sqrt{a_1}\\
0 & \displaystyle \sqrt{\frac{a_1a_2-a^2}{a_1}}
\end{matrix}
\right).
$$

\begin{Lem} \label{l61}
Suppose $A,B$ are positive real $2\times 2$ matrices, $\det A\ge 1, \det B\ge 1$, and $\det \left(\frac{A+B}{2}\right)\le 1+\delta$ for some $\delta\ge 0$. Consider $C\dd (\Lambda_A^{-1})^*B\Lambda_A^{-1}$ and write $\Lambda_C = \left(\begin{smallmatrix}
x & y\\
0 &z
\end{smallmatrix}\right)$. Then, there is $\widehat \delta$ such that
\begin{align}
&x = 1+\widehat\delta + O(\delta), \qquad z = 1-\widehat\delta + O(\delta), \qquad |y| + |\widehat\delta| = O(\sqrt{\delta}), \label{sd56} \\
&\frac{1}{4\kappa} \le x, z \le 2\sqrt{\kappa}, \qquad (2\sqrt{\kappa})^{-1} \le xz \le 2\sqrt{\kappa}, \qquad \kappa \dd  1+\delta \label{r57}.
\end{align}
Moreover, 
 \begin{equation}\label{sd_h179}
 B=\Lambda_A^*\Lambda_C^*\Lambda_C\Lambda_A, \quad \Lambda_B=\Lambda_C\Lambda_A.\end{equation}
\end{Lem}
\beginpf Identities \eqref{sd_h179} are straightforward. 
Using  $\det A \ge 1$ and $\det B \ge 1$, we obtain the estimates
\begin{align}\label{lv_1}
\det \frac{I_{2\times 2} + C}{2} \le \det A \cdot \det \frac{I_{2\times 2} + C}{2} =    \det \frac{A+\Lambda_A^*C\Lambda_A}{2} =\det \frac{A+B}{2} \le 1+\delta, \\
\det C = \frac{\det B}{\det A} \ge \frac{1}{\det A \cdot \det B} \stackrel{\rm Corollary\, \ref{sd55}}{\ge} \frac{1}{(1+\delta)^2} \geq 1-2\delta.\label{lv_2}
\end{align}
Let $C=\left(\begin{smallmatrix}u&v\\v&w\end{smallmatrix}\right)$, so 
\begin{equation}\label{lv_11}
x = \sqrt{u}, \quad y = v/\sqrt{u},\quad z = \Bigl((uw - v^2)/u\Bigr)^{1/2}.
\end{equation} 
We have 
$$
\det (I_{2\times 2}+C)=1+x^2+y^2+z^2+x^2z^2\stackrel{\eqref{lv_1}}{\leq} 4\kappa,
$$
hence 
\begin{equation}\label{lv_3}
\max(x, |y|, z, xz,\sqrt{x^2+z^2}) \le 2\sqrt{\kappa}.
\end{equation}
 Since $(xz)^2=\det C=\det B/\det A$ and $\det B\ge 1$, we also have
$$
\frac{1}{4 \kappa}\le \frac{1}{4}\det\left(\frac{A+B}{2}\right)^{-1} \le \frac{1}{\det A}\le (xz)^2,
$$ 
where the second inequality follows from Lemma \ref{sd_h51}.
Together with \eqref{lv_3}, this yields 
\begin{equation}\label{lv_4}
(2\sqrt{\kappa})^{-1}\le xz\le  2\sqrt{\kappa}.
\end{equation} Inequalities $(4\kappa)^{-1}\le x^2z^2$ and $x^2+z^2 \stackrel{\eqref{lv_3}}{\le} 4\kappa$ imply $(4\kappa)^{-1} \le 4\kappa x^2$. Therefore, $x\ge (4\kappa)^{-1}$. In a similar way, one gets $z\ge (4\kappa)^{-1}$. Thus, \ \eqref{r57} is proved. 
Note that for $\delta \ge 10^{-4}$ relations \eqref{sd56} follow from \eqref{r57} and the fact that $|y| \le 2\sqrt{\kappa}$. Now, assume that $\delta<10^{-4}$ and write
\begin{align}
uw - v^2 &= \det C \stackrel{\eqref{lv_2}}{\ge} 1 - 2\delta,\label{lv_9}\\ 
(1+u)(1+w) - v^2 &= \det (I_{2\times 2}+C) \stackrel{\eqref{lv_1}}{\le} 4(1+\delta).\notag
\end{align}
These inequalities imply 
\begin{align}
uw &\ge 1 - 2\delta, \label{lv_7}\\
u+w &\le 2(1+3\delta),\label{lv_8}
\end{align} 
so that 
\begin{equation}\label{lv_12}
(\sqrt{u} - \sqrt{w})^2 \le 2(1+3\delta) - 2\sqrt{1-2\delta} \le 2(1+3\delta) - 2(1-2\delta) \le 10\delta.
\end{equation}
Since $\delta<10^{-4}$, we have $u+w \le 4$. Hence, $\sqrt u\le 2$ and
\begin{align}
u &\ge u - \sqrt{u}(\sqrt{u}-\sqrt{w})- 2|\sqrt{u}-\sqrt{w}| = \sqrt{uw} - 2|\sqrt{u}-\sqrt{w}| 
\stackrel{\eqref{lv_7}+\eqref{lv_12}}{\ge} 1-2\delta-2\sqrt{10\delta} \label{lv_13_7},\\
&\ge 1-10\sqrt{\delta}.\notag
\end{align}
Analogous estimate holds for $w$. Moreover, we have
$$
\max(u,w) \le u + w - \min(u,w) \stackrel{\eqref{lv_8}+\eqref{lv_13_7}}{\le} 2(1+3\delta) - (1-10\sqrt{\delta}) \le 1 + 16\sqrt{\delta}.
$$
It follows that $\max(|u - 1|,|w - 1|) \le 16\sqrt{\delta}$. Since  $1 - 2\delta \stackrel{\eqref{lv_7}}{\le} \sqrt{uw} \le \frac{u+w}{2} \stackrel{\eqref{lv_8}}{\le} 1+3\delta$ and $v^2 \stackrel{\eqref{lv_9}}{\le} uw-1+2\delta$, we also have 
$$
u=1+\eps_1, \qquad w=1+\eps_2, \qquad v=\eps_3, \qquad |\eps_{1,2}|\le 16\sqrt{\delta},\qquad |\eps_1+\eps_2| \le 6\delta,
$$
and
$$
\eps_3\le \sqrt{uw-1+2\delta}\le\sqrt{(1+3\delta)^2-1+2\delta}\le 16\sqrt\delta.
$$
Relations \eqref{sd56} now follow from \eqref{lv_11} and Taylor expansion. \qed

\medskip

\noindent{\bf Remark.} In the case $\delta<10^{-4}$, the above  calculations provide explicit bounds:
$$
|u-1|\le 0.16, \quad |w-1|\le 0.16, \quad |v|\le 0.16.
$$
Thus,
\begin{equation}\label{sd_hh2}
x=\sqrt u\in (0.8,1.1), \quad z=\left(  \frac{uw-v^2}{u}\right)^{1/2}\in (0.8,1.1).
\end{equation}

\noindent{\bf Proof of Theorem \ref{t6}.} For integer $n \ge 0$, we introduce $H_n$, $\eps_n$ as follows:
\begin{equation}\label{sd_sd1}
H_n = \int_{n}^{n+1}\Hh(t)\,dt, \qquad 1 + \eps_n = \det\left(\frac{H_{n}+H_{n+1}}{2}\right).
\end{equation}
The inequality \eqref{sd_h9} from Lemma \ref{sd_mg2} and Corollary \ref{sd55} yield
\begin{equation}\label{lv_13}
\det H_n \ge 1, \quad\eps_n \ge 0
\end{equation} for all integers $n \ge 0$.
Let $\Lambda_0$ be real upper-triangular matrix (check the formula \eqref{sd_f7}) such that $H_0 = \Lambda_0^* \Lambda_0$.   Iteratively applying Lemma \ref{l61} and taking  $A = H_{n-1}$, $B = H_{n}$, $n \ge 1$, we obtain representation 
$H_n = G_n^* G_n$ for some $G_n = \Lambda_n \ldots \Lambda_0$, where $\{G_k\}$ and $\{\Lambda_k\}$, $k \ge 0$ are real matrices written coordinate-wise as
$$
\Lambda_{k+1}\dd \begin{pmatrix}
x_k & y_k\\
0 &z_k
\end{pmatrix},\quad G_{k}\dd \begin{pmatrix}
g_{1,k} & g_{2,k}\\
0 &g_{3,k}
\end{pmatrix},
$$
and satisfying
\begin{equation}\label{sa200}
\begin{aligned}
&x_k = 1+\widehat\eps_k + O(\eps_k), \qquad z_k = 1-\widehat\eps_k + O(\eps_k), \qquad |y_k| + |\widehat\eps_k| = O(\sqrt{\eps_k}),  \\
&1/(4\kappa_k) \le x_k, z_k \le 2\sqrt{\kappa_k}, \qquad (2\sqrt{\kappa_k})^{-1} \le x_k z_k \le 2\sqrt{\kappa_k}, \qquad \kappa_k \dd 1+\eps_k.
\end{aligned}
\end{equation} 
For $t \in [n, n+1), n\in \mathbb{Z}^+$, we introduce the following functions which will be used later in the proof:
\begin{align*}
Q_n(t)&= (G_n^*)^{-1}\Hh(t) G_n^{-1},\\
g(t)&= G_n+(t-n)(G_{n+1}-G_n),\\
Z(t)&= (G_{n+1}- G_n)g^{-1}(t),\\
\widehat Z(t) &= Z(t)-(\trace Z(t)/2) I_{2\times 2},\\
\nu(t) &= (\det H_0)^{-1/4} \exp\left(-\frac{1}{2}\int_{0}^{t}\trace Z(s)\,ds\right), \\  
G(t) &=  \nu(t) g(t), \\
Q(t) &=   \nu^{-2}(t)  (G_n g^{-1}(t))^*Q_n(t) (G_n g^{-1}(t)),\\ 
V(t) &=  -J\widehat Z(t),
\end{align*}
so $\Hh=G^*QG$. In the next  lemmas, we prove that it is indeed a required factorization for $\Hh$.  

\medskip

\begin{Lem}\label{l54} 
We have 
\begin{eqnarray}\label{lv_88}
Q_n(t) \ge 0,\quad \det Q_n(t) = \frac{1}{\det H_n}, \quad \int_{n}^{n+1} Q_n(t)\,dt = I_{2\times 2}, \\
\label{sd_fo20}
1\le \det H_n\le \min\Big\{(1+\eps_n)^2,4(1+\eps_n)\Big\},\\
 \max\Big\{(1+\eps_n)^{-2 }, (4(1+\eps_n))^{-1}\Big\}\le \det Q_n(t) \le 1, \label{lv_44}
\end{eqnarray}
for all $n \ge 0$, $t \in [n, n+1)$. 
\end{Lem}
\beginpf By construction, $Q_n(t) \ge 0$ and $\det Q_n(t) =1/\det H_n$. In particular, $\det Q_n$ is constant on each $[n,n+1)$. We also have 
$$
\int_n^{n+1}Q_n \, dt=(G_n^*)^{-1}\left(\int_n^{n+1} \Hh \,dt\right)G_{n}^{-1}=(G_n^*)^{-1}H_n G_{n}^{-1}=I_{2\times 2}.
$$
By \eqref{lv_13}, $\det H_n\ge 1$.  Corollary \ref{sd55} yields the bound 
$$
\sqrt{\det H_n}\le \sqrt{\det H_n \det H_{n+1}} \le \det\left(\frac{H_{n}+H_{n+1}}{2}\right)=1+\eps_n, 
$$ which gives  $\det H_n \le (1+\eps_n)^2$.
 Lemma \ref{sd_h51} provides inequality 
$$
\qquad \det H_n \le \det (H_n+H_{n+1})=4\det\left(\frac{H_{n}+H_{n+1}}{2}\right)=4(1+\eps_n),
$$
and this establishes an alternative estimate. The proof of \eqref{sd_fo20} is finished. The bounds for $\det H_n$ imply inequalities for $\det Q_n$ since $\det Q_n=1/\det H_n$.
\qed 

\medskip

\begin{Lem}\label{l56}
For every $t \ge 0$, the matrix $g(t)$ is invertible, absolutely continuous, and 
\begin{equation}\label{lv_21}
g'(t) = Z(t) g(t)
\end{equation} 
for almost every $t\in \R_+$. We also have
\begin{equation}\label{sd_f20}
Z(n+t)=\begin{pmatrix}
(x_{n}-1) u_n&   y_nu_n v_n  \\
0 & (z_{n}-1)v_n
\end{pmatrix}, \qquad 0 \le t < 1, \qquad n\ge 0,
\end{equation}
where 
\begin{equation}\label{sd_ff2}
u_{n} = 1/(1-t+tx_{n}),\quad v_n = 1/(1-t+tz_{n}).
\end{equation}
\end{Lem}
\beginpf For $t \in [0,1)$, we have 
\begin{equation}\label{lv_47}
g(t+n) = (I_{2\times 2}+t(\Lambda_{n+1}-I_{2\times 2}))G_n,\end{equation} and 
$$
\det (I_{2\times 2}+t(\Lambda_{n+1}-I_{2\times 2})) = (1-t+tx_n)(1-t+tz_n) >0, \qquad 0 \le t < 1,
$$ 
hence $g(t+n)$ is invertible and $Z(t)$ is defined correctly on $\R_+$. Direct calculation shows that $g'=Zg$ almost everywhere on $\R_+$. We also have
\begin{align}
Z(n+t)
&=(\Lambda_{n+1}-I_{2\times 2}) G_n g^{-1}(n+t) = (\Lambda_{n+1}-I_{2\times 2}) G_n (G_n+t(\Lambda_{n+1}-I_{2\times 2})G_n)^{-1}\notag\\
&=(\Lambda_{n+1}-I_{2\times 2})(I_{2\times 2}+t(\Lambda_{n+1}-I_{2\times 2}))^{-1}, \label{sa204}
\end{align}
which yields  \eqref{sd_f20}. \qed

\medskip

\noindent {\bf Remark.} This lemma allows us to write
\begin{equation}\label{sd_f70}
\det g(t)=\det g(0)\cdot \exp \left(\int_0^t \trace Zd\tau\right)
\end{equation}
and we can use definition of $\nu$ to get
\begin{equation}\label{sd_f21}
\nu(t)=\frac{\sqrt{\det g(0)}}{\sqrt{\det g(t)}}\cdot 
\frac{1}{(\det H_0)^{1/4}}=(\det g(t))^{-1/2},
\end{equation}
since $\det g(0)=(\det H_0)^{1/2}$. We notice here that $\det g(t)>0$ for all $t$ thus we can take a square root in the formula above.
\begin{Lem}
For $n \ge 0$, $t \in [n,n+1)$, we have 
\begin{equation}\label{lv_32}
\trace Z(t) = O(\eps_n).
\end{equation}
\end{Lem}
\beginpf For $t\in [0,1)$, we have 
$$
|\trace Z(n+t)| \stackrel{\eqref{sd_f20}}{=} |(x_n-1)u_n + (z_n - 1)v_n| =\left|\frac{x_n-1}{1-t+tx_n}+\frac{z_n-1}{1-t+tz_n}\right|.
$$
If $\eps_n<10^{-4}$, we use  \eqref{sd_hh2}, \eqref{sa200}, and Taylor expansion to get
$$
\frac{x_n-1}{1+t(x_n-1)}+\frac{z_n-1}{1+t(z_n-1)}=x_n+z_n-2+O(|x_n-1|^2+|z_n-1|^2)=O(\eps_n).
$$
If $\eps_n>10^{-4}$, we recall \eqref{sa200} again and  write
\begin{equation}\label{sd_ff1}
\left|\frac{x_n-1}{1+t(x_n-1)}+\frac{z_n-1}{1+t(z_n-1)}\right|\le |x_n-1|+|x_n^{-1}-1|+|z_n-1|+|z_n^{-1}-1|=O(\eps_n),
\end{equation}
where we, first, used the fact that a linear function in $t\in [0,1]$ achieves its minimum at an endpoint of the segment $[0,1]$ and, second, combined all four possible values in the sum in right hand side. \qed 

\medskip

\begin{Lem}\label{l57}
The matrix-function $V$ is symmetric and has real entries  for all $t \in \R_+$. Moreover, there exist $V_{1, }$, $V_{2}$ such that $V = V_{1} + V_{2}$, and $\|V_{1}\|_{L^1} \lesssim
\widetilde \K({\Hh})$ and
$\|V_{2}\|_{L^2}^{2} \lesssim \widetilde \K({\Hh})$.
\end{Lem}
\beginpf Indeed, for all $t \in [0, 1)$, we have 
\begin{equation}\label{lv_101}
V(n+t) = -J\widehat Z(n+t) = \frac{1}{2}
\begin{pmatrix}
0 &   (z_n - 1)v_n-(x_n-1)u_n  \\
(z_n - 1)v_n-(x_n-1)u_n & -2y_nu_n v_n
\end{pmatrix} = V^*(n+t).
\end{equation}
It follows from \eqref{sa200} that $|x_n -1|+ |z_n -1| \lesssim \max\{\sqrt{\eps_n},\eps_n\}$ and $|y_n| \lesssim \max\{\sqrt{\eps_n},\eps_n\}$ for every $n \ge 0$. In particular,  there exists $\widetilde\eps>0$ such that the estimate
$\eps_n<\widetilde\eps$ implies 
$$
|x_n - 1|+ |z_n - 1| \le 1/2, \quad u_n\lesssim 1,\quad v_n \lesssim 1.
$$
Therefore, if we define $\zeta_n = \chi_{\{n:\eps_n<\widetilde\eps\}}$ and
$d_{2, n} = \zeta_n y_nu_n v_n$, then $|d_{2,n}|\lesssim \sqrt{\eps_n}$.
Let  $d_{1,n} =  (1-\zeta_n)y_nu_nv_n=y_nu_n v_n - d_{2,n}$ and,  finally, define $V_{1}$ and $V_2$ on each $[n, n+1)$ by (recall that $\trace Z=(x_n-1)u_n+(z_n-1)v_n)$, see \eqref{sd_f20}) 
\begin{align*}
V_1 &= 
\begin{pmatrix}
0 &   \trace Z/2-(x_n-1)(1-\zeta_n)u_n  \\
\trace Z/2-(x_n-1)(1-\zeta_n)u_n & -d_{1,n}
\end{pmatrix},\\
V_2 &= 
-\begin{pmatrix}
0 &   (x_n-1)\zeta_n u_n  \\
(x_n-1)\zeta_n u_n & d_{2,n}
\end{pmatrix},
\end{align*}
so the identity $V=V_1+V_2$ follows from \eqref{lv_101}.
By \eqref{lv_32}, we have 
\begin{equation}\label{lv_77}
\|\trace Z\|_{L^1} \lesssim \widetilde\K({\Hh}).
\end{equation}
Similarly to \eqref{sd_ff1}, we get
\begin{equation}\label{lv_88n}
|(x_n-1)u_n(1-\zeta_n)| \stackrel{\eqref{sa200}}{\lesssim} \eps_n.
\end{equation}
We claim that
\begin{equation}\label{lv_110}
|d_{1,n}| {\lesssim } \,\eps_{n} 
\end{equation}
for every $t\in [0,1]$. Indeed, we have 
$$
d_{1,n}=
 \left\{\begin{array}{cc} y_nu_n v_n, & n:  \eps_n > \widetilde\eps,\\
0, &  n: \eps_n \le \widetilde\eps\end{array}\right.,\quad y_nu_nv_n\stackrel{\eqref{sd_ff2}}{=}\frac{y_n}{(t(x_n-1)+1)(t(z_n-1)+1)}.
$$
At the endpoints of $[0,1]$ the quadratic polynomial $P(t)\dd (t(x_n-1)+1)(t(z_n-1)+1)$ takes values~$1$ or $x_nz_n$.  It is also positive on $[0,1]$. We can use \eqref{sa200} to write a bound $x_nz_n\gtrsim \eps_n^{-1/2}$. If the first coefficient satisfies $(x_n-1)(z_n-1)\le 0$, $P$ reaches minimum over $[0,1]$ at an endpoint and we are done because $|y_n|\lesssim \sqrt{\eps_n}$ and $|y_n/(x_nz_n)|\lesssim \eps_n$. Otherwise, consider, e.g., the case $x_n\in (0,1),z_n\in (1,\infty)$.
The point of minimum  of $P$ over $\R$ is given by
$\displaystyle 
t_*=-\frac{x_n+z_n-2}{2(x_n-1)(z_n-1)}.
$
If $x_n+z_n-2\ge 0$, then $t_*\le 0$ and $P$ again reaches minimum over $[0,1]$ at zero, an endpoint of $[0,1]$. If, however, $x_n+z_n-2< 0$, we get $z_n<2$ and $x_n\ge (x_nz_n)/2\gtrsim \eps_n^{-1/2}$ so we can write
$$
(t(x_n-1)+1)(t(z_n-1)+1)\ge \min\{ 1,x_n,z_n,x_nz_n\}=x_n\gtrsim \eps_n^{-1/2}.
$$
Thus, we have \eqref{lv_110} in all cases. Summarizing, we get
$$
\int_{0}^{+\infty}\|V_1(t)\|\,dt \stackrel{\eqref{lv_77}+\eqref{lv_88n}+\eqref{lv_110}}{\lesssim}  \widetilde \K({\Hh}). 
$$ 
Since  $|d_{2,n}| \lesssim \sqrt{\eps_n}$ 
and
$
\zeta_n|x_n-1|u_n \stackrel{\eqref{sa200}}{\lesssim} \sqrt{\eps_n},
$
we have
$$
\int_{0}^{+\infty}\|V_2(t)\|^2\,dt \lesssim \sum_{n \ge 0}\eps_n= \widetilde \K({\Hh}). 
$$ 
The lemma follows. \qed

\begin{Lem} \label{sd_l58} We have $\det G=\det Q=1$ and
\begin{equation}
\|\trace Q-2\|_{L^1(\mathbb{R}^+)}\lesssim \widetilde \K(\Hh).
\end{equation}
\end{Lem}
\beginpf
Notice that $\det G =\nu^2(t)\cdot \det g \stackrel{\eqref{sd_f21}}{=}  1$. Since  $\Hh=G^*QG$ and $\det \Hh=1$, we also have $\det Q=1$.
Recall  that $Q>0$ and $\trace Q(t) - 2 \ge 2\sqrt{\det Q(t)} - 2 = 0$. So, we only need  an estimate for $\trace Q(t)$ from above. 
For each $n\in \mathbb{Z}^+$, we consider 
$$
\int_{n}^{n+1} (\trace Q-2)\, dt
$$
and handle separately the cases of small and large $\eps_n$.

\medskip

\noindent {\bf Case 1.} Assume that $\eps_n<1$. Define $T_n: t \mapsto G_n g^{-1}(n+t)$ on $[0, 1)$. Then, for $t \in [0, 1)$, we use \eqref{sd_f21} to write
\begin{eqnarray}\label{lv_200}
\trace Q(t) = (\det g)\cdot  \trace\bigl(T^*_{n}(t) Q_n(t) T_n(t)\bigr) = (\det g)\cdot \trace\bigl(T_n(t)T^*_{n}(t)Q_n(t)\bigr).\notag
\end{eqnarray}
From \eqref{sd_f70}, we get
$$
\frac{\det g(n+t)}{\det g(n)}=\exp\left(\int_n^t \trace Zd\tau\right).
$$
Since $g(n)=G_n, \det G_n=(\det H_n)^{1/2}$, we recall \eqref{sd_fo20} and \eqref{lv_32} to write
\begin{equation}\label{sd_ff89}
\det g(n+t)=(1+O(\eps_n))\exp(O(\eps_n))=1+O(\eps_n).
\end{equation}
For $t \in [0, 1)$, we have $T_n(n+t) \stackrel{\eqref{lv_47}}{=} G_n (G_n+t(\Lambda_{n+1}-I_{2\times 2})G_n)^{-1} = (I_{2\times 2} + t(\Lambda_{n+1}-I_{2\times 2}))^{-1}$, that is,
$$
T_n(n+t) 
= \begin{pmatrix}u_n& -ty_{n}u_nv_{n} \\ 0 &v_n\end{pmatrix}, \qquad \begin{cases}\det (T_nT_n^*) = u_n^2 v_n^2, \\ 
\trace (T_nT_n^*)= u_n^2 + v_n^2 + (t y_{n} u_n v_n)^2.\end{cases}
$$
Recalling the formulas \eqref{sd_ff2} for $u_n,v_n$, we get
$$
u_n^2 v_n^2 \stackrel{\eqref{sa200}}{=} 1 + O(\eps_n), \qquad u_n^2 + v_n^2 + (t y_{n} u_n v_n)^2 \stackrel{\eqref{sa200}}{=} 2 + O(\eps_n).
$$
We remind that $\lambda_{\min}(A)$ and $\lambda_{\max}(A)$ denote smallest and largest eigenvalues of self-adjoint matrix $A$, respectively.
Then, the formulas for $u_n$, $v_n$ show that 
$$
\lambda_{\min}(T_nT_n^*)\lambda_{\max}(T_nT_n^*)=1+O(\eps_n),\quad 
\lambda^2_{\min}(T_nT_n^*)+\lambda^2_{\max}(T_nT_n^*)
=2+O(\eps_n).
$$
Hence,
$$
 (\lambda_{\min}(T_nT_n^*)+ \lambda_{\max}(T_nT_n^*)   )^2=4+O(\eps_n),\quad (\lambda_{\min}(T_nT_n^*)- \lambda_{\max}(T_nT_n^*)  )^2=O(\eps_n),
$$
and
$$
\lambda_{\max}(T_nT_n^*) = 1 + \kappa_n(t) + O(\eps_n), \quad \lambda_{\min}(T_nT_n^*)= 1 - \kappa_n(t) + O(\eps_n),
$$ 
for some function $\kappa_n \ge 0$ satisfying $\kappa_n= O(\sqrt{\eps_n})$ on $[0, 1)$. Then,
$$
\int_{n}^{n+1}(\lambda_{\max}(Q_n)+ \lambda_{\min}(Q_n))\,dt = \int_{n}^{n+1} \trace Q_n(t)\, dt = \trace \int_{n}^{n+1} Q_n(t)\,dt \stackrel{\eqref{lv_88}}{=} 2,   
$$
so Cauchy-Schwarz inequality implies
$$
\left(\int_n^{n+1} ( \lambda_{\max}(Q_n)  -  \lambda_{\min}(Q_n)   )\,dt\right)^2 \le 4\int_n^{n+1} \left(\sqrt{ \lambda_{\max}(Q_n)} - \sqrt{  \lambda_{\min}(Q_n)  }\right)^2\,dt.
$$
The right hand side of the above formula equals 
$$
4 \left(2 - 2\int_{n}^{n+1}\sqrt{\det Q_n}\,dt\right) = O(\eps_n),
$$
as follows from Lemma \ref{l54}. Using von Neumann inequality for the trace of a product of two matrices \cite{mir}, we obtain
\begin{align*}
\int_{n}^{n+1}\trace (T_n(t)T_{n}^*(t)Q_n(t))\,dt 
&\le \int_{n}^{n+1}( \lambda_{\max}(T_nT_n^*)  \lambda_{\max}(Q_n) +  \lambda_{\min}(T_nT_n^*)  \lambda_{\min}(Q_n))\,dt,  \\
&\le 2 + \sup_{[n, n+1)}\kappa_n(t)\cdot\int_{n}^{n+1}\Bigl(\lambda_{\max}(Q_n)- \lambda_{\min}(Q_n)\Bigr)\,dt +  O(\eps_n), \\
&\le 2 +  O(\eps_n).
\end{align*}
We now use \eqref{sd_ff89} to obtain $\displaystyle \int_{n}^{n+1} (\trace Q-2) \, dt\lesssim \eps_n.$

\noindent {\bf Case 2.} Assume that $\eps_n\ge 1$. We only need to show that
$\displaystyle 
\int_{n}^{n+1} \trace Q\, dt\lesssim \varepsilon_n
$
since that would imply the bound
$$
\sum_{n:\eps_n\ge 1}\int_{n}^{n+1} (\trace Q-2)\,dt\lesssim \sum_{n\ge 0}\eps_n .
$$
We can write $Q=(\det g) T_n^*Q_nT_n$. Then, notice that
$$
\trace Q=(\det g)\cdot \trace (T_nT_n^*Q_n)\le (\det g)\cdot \lambda_{\max}(T_nT_n^*)\trace Q_n,
$$
again by von Neumann inequality for the trace. Introducing 
$$
Y_n(n+t) =  T_n^{-1}(n+t)=((1-t)I_{2\times 2}+tG_{n+1}G_n^{-1}),
$$ 
we can write
$$
\det g=\det G_n\cdot \det Y_n=\frac{\det G_n}{ \sqrt{\det (T_nT_n^*)}}=\frac{\det G_n}{ \sqrt{\lambda_{\min}(T_nT_n^*)\lambda_{\max}(T_nT_n^*)}},
$$
so  
\begin{align*}
(\det g)\cdot \lambda_{\max}(T_nT_n^*)\trace Q_n
&=(\det H_n)^{1/2}\sqrt{\frac{\lambda_{\max}(T_nT_n^*)}{\lambda_{\min}(T_nT_n^*)
}}\trace Q_n,\\
&=(\det H_n)^{1/2}\sqrt{\frac{\lambda_{\max}(Y_n^*Y_n)}{
\lambda_{\min}(Y_n^*Y_n)}}\trace Q_n.
\end{align*}
Since $\int_n^{n+1}Q_ndt=I_{2\times 2}$, we only need to show that
\begin{equation}\label{sd_f2}
(\det H_n)^{1/2}\sqrt{\frac{\lambda_{\max}(Y_n^*Y_n)}{
\lambda_{\min}(Y_n^*Y_n)}}\lesssim \varepsilon_n.
\end{equation}
Let us apply Lemma \ref{sd_f1}. We take $A=H_n$,$B=(G_n^*(1-t)+G_{n+1}^*t)(G_n(1-t)+G_{n+1}t), \Omega=Y_n^*Y_n$ and notice that $Y_n^*Y_n=(G_n^*)^{-1}BG_n^{-1}=(\Lambda_A^*)^{-1}B\Lambda_A^{-1}$ and $\alpha=\det H_n$. 
Let us estimate quantities $\beta,\gamma$ from Lemma \ref{sd_f1}. We get
$$
\det (G_n(1-t)+G_{n+1}t)=(g_{1,n}(1-t)+g_{1,n+1}t)(g_{3,n}(1-t)+g_{3,n+1}t)\ge (1-t)^2g_{1,n}g_{3,n}+t^2g_{1,n+1}g_{3,n+1},
$$
since the diagonal elements of $G_n$ and $G_{n+1}$ are positive by definition. 
Moreover $g_{1,n}g_{3,n}=\det G_n=(\det H_n)^{1/2}\ge 1$ and 
$g_{1,n+1}g_{3,n+1}=\det G_{n+1}=(\det H_{n+1})^{1/2}\ge 1$ so $\beta=\det B\ge 1/4$.  For the quadratic form, we have an estimate
$$
\langle A\xi,\xi\rangle+\langle B\xi,\xi\rangle\lesssim \langle H_n\xi,\xi\rangle+\langle H_{n+1}\xi,\xi\rangle+\Re\,\langle G_n^*\xi,G_{n+1}\xi\rangle \lesssim \langle H_n\xi,\xi\rangle+\langle H_{n+1}\xi,\xi\rangle, \quad \xi\in \C^2\,,
$$
in which we applied Cauchy-Schwarz inequality in the second bound.
So,
$
A+B\le C(H_n+H_{n+1})
$. Thus, $\gamma =\det (A+B)\lesssim \det (H_n+H_{n+1})\stackrel{\eqref{sd_sd1}}{\lesssim} \varepsilon_n$ since $\eps_n\ge 1$. Now, Lemma \ref{sd_f1} gives \eqref{sd_f2}.
\qed

\medskip

We are ready to complete the proof of Theorem \ref{t6}. From the definition of $G$, $Q$, we see that $\Hh = G^* Q G$ on $\R_+$ and we already established that $\det G=\det Q=1$. Moreover, 
$$
G'=\nu'g+\nu g'=\nu'g+\nu Zg=\nu\left(Z-\frac{1}{2}\trace Z\cdot I_{2\times 2}\right)g=JVG,\quad V=V^*.
$$ 
Lemma \ref{l57} and Lemma \ref{sd_l58} provide the necessary bounds for $V$ and $Q$ and that finishes the proof.
 \qed

\medskip

\section{Appendix}\label{s10}
In this Appendix, we collect some auxiliary statements used in the main text. Most of them are well-known but we give proofs for the reader's convenience. 

\begin{Lem}\label{sd_h7}
If $A\in {\rm SL}(2,\R)$, then $A^*JA=J, \quad A^{-1}=-JA^*J, \quad JAJ^{*}=(A^*)^{-1}$.
\end{Lem}
\beginpf
The proof is  a straightforward calculation.
\qed\bigskip

We recall that $\Lambda_A$ denotes the upper-triangular matrix providing factorization of matrix $A$ introduced in~\eqref{sd_f7}. 

\begin{Lem}\label{sd_f1} Suppose $A$ and $B$ are real positive symmetric $2\times 2$ matrices, $\det A= \alpha$, $\det B=\beta$, and
$\det (A+B)=\gamma$. If $\Omega = (\Lambda_A^*)^{-1}B\Lambda_A^{-1}$, then
$$
\alpha\frac{\lambda_{\max}(\Omega)}{\lambda_{\min}(\Omega)}\le \frac{\gamma^2}{\beta}.
$$
\end{Lem}
\beginpf Denote $x=\lambda_{\min}(\Omega), y=\lambda_{\max}(\Omega)$, $t=y/x$. Clearly, $x,y>0, t\ge 1$. Then, we have
$$
\beta=\det B=\det (\Lambda^*_A\Omega \Lambda_A)=\alpha x^2t, \quad \det (A+B)=\alpha\det (I_{2\times 2}+\Omega)=\alpha(1+x)(1+xt)=\gamma.
$$
Thus, $tx\le \gamma/\alpha$ and $t\sqrt{\beta/(\alpha t)}\le \gamma/\alpha$ so $\alpha t\le \gamma^2/\beta$. \qed

\begin{Lem}\label{sd_ak3}
Suppose $\Omega$ is a matrix-function defined on $\R_+$ and integrable over any finite interval. Denote the largest eigenvalue of $\Omega(t)+\Omega^*(t)$ by $\Lambda(t)$. If $X(t)$ is absolutely continuous vector-function that solves
$
X'=\Omega X,
$
then $\|X(t)\|\le \|X(0)\|\exp\left(\displaystyle \frac{1}{2}\int_0^t \Lambda(\tau)d\tau\right)$.
\end{Lem}
\beginpf
If $\Psi = \|X\|^2$, then
$$
\Psi'=\langle \Omega X,X\rangle+\langle X,\Omega X\rangle=\langle (\Omega+\Omega^*) X,X\rangle\le \Lambda\|X\|^2=\Lambda\Psi
$$
and we get statement of the lemma.
\qed

\begin{Lem}\label{sd_h51}If $A$ and $D$ are two real non-negative $2\times 2$ matrices and $A\le D$, then
$$
\det A\leq \det D.
$$
Moreover, if we have equality above and $\det A>0$, then $A=D$. 
If we have equality and $\det A=0$, then either $A=0$ or there is $\mu\in [1,\infty)$ such that $D=\mu A$.
\end{Lem}
\beginpf If $\det D>0$, the statement becomes trivial after we notice that $A\le D$ is equivalent to $D^{-1/2}AD^{-1/2}\le I_{2\times 2}$. If $\det A=\det D=0$ and $A\le D$, then 
$$
U^{-1}AU\le \left(
\begin{array}{cc}
\lambda_1 &0\\
0& 0
\end{array}
\right),
$$
where unitary $U$ diagonalizes $D$ and $\lambda_1$ is an eigenvalue of $D$. This implies that $D=\mu A$ with some $\mu\in [1,\infty)$. \qed

\medskip

The following bound is known as Minkowski estimate for determinants (e.g., \cite{marcus}, p.115).

\begin{Lem}\label{lr}If $A$ and $B$ are two non-negative real $2\times 2$ matrices, then
$$
\det (A+B)\ge (\sqrt{\det A}+\sqrt{\det B})^2,
$$
and equality holds if and only if one of the following conditions holds:
\begin{itemize}
\item $\det B>0$ and $A=\mu B$ with some $\mu\in [0,\infty)$, 
\item $\det A>0$ and $B=\mu A$ with some $\mu\in [0,\infty)$,
\item $\rank B+\rank A\le 1$,
\item $\rank B=\rank A=1$ and there is $\kappa\in (0,\infty)$ such that $A=\kappa B$.
\end{itemize}
\end{Lem}
\beginpf
If $\det B=0$ or $\det A=0$, the proof follows from the previous lemma. Otherwise, we can always reduce the setup to the case when $B=I$ by dividing the both sides by $\det B$. 
If
$$
A=\begin{pmatrix}
a_1&a\\
a&a_2
\end{pmatrix}, 
$$
we only need to check that
$$
(1+a_1)(1+a_2)-a^2\ge 1+a_1a_2-a^2+2\sqrt{a_1a_2-a^2},
$$
which is equivalent to 
$$
(a_1-a_2)^2\ge -4a^2.
$$
Equality holds if and only if $A=\kappa I, \kappa\in [0,\infty)$. \qed

\medskip

We immediately get the following corollary.
\begin{Cor}\label{sd55} Suppose $A,B$ are two real non-negative $2\times 2$ matrices, then
\begin{eqnarray}\notag
\det\frac{A+B}{2}\geq \sqrt{\det A\det B},\quad
\det (A+B)\ge \det A+\det B.
\end{eqnarray}
\end{Cor}
\begin{Lem}\label{sd_mg2}
Let $H$ be real and non-negative $2\times 2$ matrix function on $[a,b]$, $\det H(t)=1$ for almost every $t\in [a,b]$, and $H\in L^1(a,b)$. Then, we have
\begin{eqnarray}\label{sd_h9}
\det \int_a^b H(t)\, dt \ge (b-a)^2.
\end{eqnarray}
Moreover, equality holds if and only if $H$ is constant almost everywhere on $[a,b]$.
\end{Lem}
\beginpf  By a change of variables, we can reduce the problem to the case when $a=0,b=1$.  We have  
$$
\frac{1}{\sqrt{\det A}}=\frac{1}{\pi}\int_{\mathbb{R}^2} e^{- \langle Ax,x \rangle}\, dx
$$
for every real matrix $A$ with nonzero determinant. Take $A = \displaystyle \int_0^1 H\, dt$. By Jensen's inequality, we have
$$
\frac{1}{\pi}\int_{\mathbb{R}^2} e^{- \langle Ax,x \rangle}\, dx \le \int_0^1  \left(\frac{1}{\pi}\int_{\mathbb{R}^2} e^{-\langle H(t)x,x \rangle}dx\right)dt=\int_0^1 \frac{1}{\sqrt{\det H(t)}}dt=1.
$$
If equality holds in \eqref{sd_h9}, then $\langle H(t)x,x\rangle$ is constant in $t$ for every $x\in \R^2$. We can call this constant $C_x$. By polarization identity,  $C_{x+y}-C_{x-y}=4\langle H(t)x,y\rangle$ is constant in $t$ for every $x,y$. Taking $x=e_j$, $y=e_l$ for $j,l=1,2$, where $\{e_n\}$ is standard basis in $\R^2$, we see that elements of $H$ are constants in $t$.
\qed
\begin{Lem}\label{sd_h1}
Let $H(t)$ be real and non-negative $2\times 2$ matrix function and $H\in L^1(a,b)$. Then, 
\begin{equation}\label{sd_h10}
\det \int_a^b H(t)\, dt\ge \left(\int_a^b \sqrt{\det H(t)}\,dt\right)^2.
\end{equation}
Assuming that $\trace H>0$ almost everywhere on $[a,b]$, we have equality in \eqref{sd_h10} if and only if $H$ is equivalent to a non-negative constant matrix $\displaystyle \int_a^b H(t)\, dt$.
\end{Lem}
\beginpf We  can assume that $a=0$ and $b=1$. Let us first do the proof assuming that there is $\delta>0$ such that 
\begin{equation}\label{122}
\det H(t)> \delta, \quad t\in [0,1].
\end{equation}
Consider increasing function
$$
\upsilon(t)=\int_0^t \sqrt {\det H(\tau)}\,d\tau, \quad t\in [0,1]
$$
and let $\eta$ define its inverse function so that
\begin{equation}\label{lv_111n}
\eta(0)=0, \quad \eta'(\upsilon)=\frac{1}{\upsilon'(\eta(\upsilon))}=\frac{1}{\sqrt{\det H(\eta(\upsilon))}}.
\end{equation}
We write
$$
\det \int_0^1 Hdt=\det \int_{0}^{\upsilon(1)}H(\eta(\tau))\eta'(\tau)\,d\tau.
$$
Formula \eqref{lv_111n} makes sure that the matrix under the last integral has unit determinant and the previous lemma gives
$$
\det \int_0^1 H(t)\, dt\ge \upsilon^2(1).
$$
On the other hand, 
$$
\int_0^1 \sqrt{\det H(t)}\,dt=\int_0^{\upsilon(1)} \sqrt{\det H(\eta(\tau))}\cdot \eta'(\tau)\,d\tau=\upsilon(1),
$$
since the integrand is equal to $1$. Now that we have  proved the lemma under assumption \eqref{122}, we can use the standard approximation argument (e.g., by considering $H+\delta I_{2\times 2}$ with $\delta>0$ and then sending $\delta\to 0$), to show \eqref{sd_h10} in full generality. 

\medskip 

Next, assume that $a=0$, $b=1$ and that we have equality in \eqref{sd_h10}. Then, 
for every $c\in (0,1)$, we get
\begin{flalign*}
A_c\dd\int_0^c Hdt, \quad B_c\dd\int_c^1Hdt,\\
\det (A_c+B_c)=\left(\int_0^c \sqrt{\det H(t)}\,dt+\int_c^1 \sqrt{\det H(t)}\,dt\right)^2\stackrel{\eqref{sd_h10}}{\le}& \left( \sqrt{\det A_c} +\sqrt{\det B_c}\right)^2.
\end{flalign*}
Lemma \ref{lr} provides us with an opposite bound so we actually have equality in the estimate above. 
Notice that $A_c\neq 0$ and $B_c\neq 0$ by our assumptions on the trace.  If $\det A_{c^*}=0$ for some $c^*$, then $\det B_{c^*}=0$ by Lemma \ref{lr}. Moreover, if $\det A_{c_2}=0$, then $\det A_{c_1}=0$ for all $c_1<c_2$ since $A_{c_1}\le A_{c_2}$. Similarly, if $\det B_{c_1}=0$, we get $\det B_{c_2}=0$. Thus, if $\det A_{c^*}=0$ for some $c^*$, then $\det A_c=0$ for all $c$ and, by continuity, 
$$
\det\int_0^1 H(t)\,dt=0.
$$
Then, by Lemma \ref{sd_h51}, we get
\begin{equation}\label{sd_vi2}
\alpha(c)\int_0^1 H(t)\,dt=\int_0^c H(t)\, dt,
\end{equation}
with $\alpha(c)\in (0,1)$. Taking trace of both sides, we get
$$
\alpha(c)\int_0^1 \trace H(t)\,dt=\int_0^c\trace H(t)\,dt.
$$
Therefore,
$
\displaystyle \alpha(c)=\left(\int_0^c\trace Hdt\right)\left(\int_0^1\trace Hdt\right)^{-1}$ and differentiation of \eqref{sd_vi2} in $c$ gives $$\quad H(c)=\left(\frac{\trace H(c)}{\int_0^1\trace Hdt}\right)\int_0^1 Hdt.
$$
So, $H$ is equivalent to $\displaystyle \int_0^1Hdt$.

\medskip

Let us suppose now that $\det A_c>0$ for all $c$. By Lemma \ref{lr}, there is a positive  function $\nu(c)$ such that
$
\nu(c)A_c=B_c$ and so $$(1+\nu(c))\int_0^c H(t)\,dt=\int_0^1 H(t)\,dt.
$$
In a similar way, we get
$$
H(c)=\left(\frac{\trace H(c)}{\int_0^1\trace H(t)\,dt}\right)\int_0^1 H(t)\,dt,
$$
and $H$ is equivalent to $\displaystyle \int_0^1 H(t)\,dt$. \qed

\begin{Lem}\label{sd_h13}Suppose $m\in \cal{N}(\C_+)$ and $m\ne 0$. Then, 
$$
\log|m(i)|=\frac{1}{\pi}\int_\R \frac{\log |m(x)|}{1+x^2}\,dx.
$$
\end{Lem}
\noindent We skip the proof of this well-known fact, which is  based on mean value formula for harmonic functions.

\bibliographystyle{plain} 
\bibliography{bibfile}
\end{document}